\documentclass[a4paper,11pt]{article}

\usepackage{subfigure}

\usepackage{latexsym}
\usepackage[english]{babel}
\usepackage{amsmath}
\usepackage{amssymb}
\usepackage{pifont}
\usepackage{float}
\usepackage{verbatim}
\usepackage{amsfonts,amsthm,mathrsfs,wasysym}%FRANCESCA

\usepackage[dvips]{graphicx}%FRANCESCA
\usepackage{pstricks,a4}%FRA

\usepackage{lmodern}%FRA
\usepackage[T1]{fontenc}%FRA 

\usepackage{sectsty}%FRA

\newtheorem{TEO}{Theorem}[section]
\newtheorem{PROP}[TEO]{Proposition}
\newtheorem{COR}[TEO]{Corollary}
\newtheorem{LEMMA}[TEO]{Lemma}
\newtheorem{DEF}[TEO]{Definition}

\newtheorem{REM}[TEO]{Remark}

\def\black{\ {\hbox{\vrule width 4pt height 4pt depth
0pt}}}
\def\Proof{{\sl Proof.}\quad}
\def\Proof{{\sl Proof.}\quad}

\newcommand{\fine}{~\hfill~$\black$\newline}

\def\s{\sigma}
\def\o{\omega}
\def\b{\beta}
\def\g{\gamma}

\def\ra{\rightarrow}

\def\gjm{-\gamma j\left[\frac{1}{N^{1/4}(\mathrm{ch}(\beta)+\mathrm{ch}(\gamma h))}\left(\bar{x} + \frac{\bar{v}}{2\mathrm{sh}(\gamma h)(\mathrm{ch}(\beta)+2\mathrm{ch}(\gamma h))}\right) + kh\right]}
\def\gjmm{-\gamma j\left[\frac{1}{N^{1/4}(\mathrm{ch}(\beta)+\mathrm{ch}(\gamma h))}\left(\widetilde{x}_N(t)+ \frac{\widetilde{v}_N(t)}{2\mathrm{sh}(\gamma h)(\mathrm{ch}(\beta)+2\mathrm{ch}(\gamma h))}\right) + kh\right]}

\def\aa{\mathrm{ch}(\b)+\mathrm{ch}(\g h)}
\def\bb{\mathrm{ch}(\b)+2\mathrm{ch}(\g h)}

\begin{document}

\title{A simple mean field model for social interactions: dynamics, 
fluctuations, criticality}
\author{Francesca Collet, Paolo Dai Pra, Elena Sartori
\\[0.2cm]Dipartimento di Matematica Pura ed Applicata, University of 
Padova \\
63, Via Trieste;  I - 35121 - Padova, Italy\thanks{e-mails:
fcollet@math.unipd.it, daipra@math.unipd.it,
esartori@math.unipd.it}}

\date{}

\maketitle

\begin{abstract}
\noindent We study the dynamics of a spin-flip model with a mean field 
interaction. The system is non reversible, spacially inhomogeneous, and 
it is designed to model social interactions. We obtain the limiting 
behavior of the empirical averages in the limit of infinitely many 
interacting individuals, and show that phase transition occurs. Then, 
after having obtained the dynamics of normal fluctuations around this 
limit, we analyze long time fluctuations for critical values of the 
parameters. We show that random inhomogeneities produce critical 
fluctuations at a shorter time scale compared to the homogeneous system.

\vspace{0.3cm}


\end{abstract}

\section{Introduction}

The mathematical description  of complex social systems has been largely 
inspired by modeling of physical systems, in particular by Statistical 
Physics (see e.g. \cite{MaSt99, Sta06}). Social interactions have, 
however, their own peculiar features.
\begin{itemize}
\item[$\RHD$]
In many applications, the same information is potentially available to 
all individuals; thus geometrical constraints in the interaction are not 
justified. Unlike in physical systems, interactions of {\em mean field} 
type may provide accurate descriptions of real behaviors.
\item[$\RHD$]
In physical  systems, interactions are often coded in an energy 
function, the {\em Hamiltonian}; the associated Gibbs distributions 
describe the equilibrium behavior of the system in the thermodynamic 
limit. When a stochastic dynamic model is desirable, for instance for 
particles in a heat bath, ``natural'' dynamics are obtained by adding 
stochastic perturbations to the Hamiltonian dynamics, preserving the 
reversibility with respect to the Gibbs distribution. In social systems, 
interactions may be more naturally given in a dynamic framework. For 
instance, at each time step an agent in a market changes his own 
``state'' (e.g. the amount of money invested in a specific item) to 
maximize his own {\em utility}, that depends on the global state of the 
business network he belongs to, and on some randomness. The associated 
stochastic evolutions are not necessarily reversible with respect to 
their equilibrium distribution.
\end{itemize}

\noindent
The mean field assumption could be inappropriate in many context, where individuals tend to conform their choices to those of a small number of appropriately defined neighbors. These situations are better described by models with {\em local interactions}, such as the {\em Voter Model} (see \cite{Lig99}); for applications of this model to social science see e.g. \cite{GiWe04}, \cite{GiWe06}, where the graph of the interactions is regular, and \cite{SoRe05}, \cite{SuEgSM05} for the case of interactions in an heterogeneous  graph. There are however many situations where the mean field assumption is reasonable. This is the case, for example,  when one models the economic behavior of agents sharing the same information see e.g. \cite{Blu03}, \cite{BlDu03}, \cite{BrDu01}. More recent developments along these lines are found in \cite{BaTo09}, \cite{DPRST09}, \cite{DaPTo}, \cite{FrBa08}. In these models each agent aims at optimizing an utility function, whose dependence on the states of other agents is invariant under permutation (i.e. mean field). Moreover, agents update their state sequentially; in other words, simultaneous updating is not allowed. As observed in \cite{CoLo03}, the purpose in this approach ``is not to model strategic interactions but collective behavior of non-strategic nature''. When simultaneous updating is allowed, non-cooperative behavior of agents needs to be considered; a dynamic game-theoretic framework in the context of mean field interaction has been recently developed in \cite{LaLi07}, while applications to economics are considered in \cite{LaSaTu09}.

A typical feature of dynamic mean field models is that they exhibit {\em phase transition}: in  the limit of infinitely many individuals, different initial conditions may be attracted to different equilibria. In a social setting, different equilibria may correspond to different degrees of polarization of opinions, better or worse financial state in a network of interacting agents and so on.
The model we study in this paper, that generalizes the one introduced in 
\cite{DPRST09}, should be seen as a prototype model for social systems 
having the following features.
\begin{enumerate}
\item
The dynamics concern the states of a large number, $N$, of interacting 
individuals.
\item
Let $\o_i$ be the state of the $i$-th individual. Other individuals {\em 
perceive} this state subject to a random perturbation. We denote by 
$\s_i$ the {\em perceived state}.
\item
Each individual changes his state by ``adapting'' to the perceived 
state of all the others. In other words, also considering the mean field 
assumption, the rate at which $\o_i$ changes depends on the empirical 
mean of the $\s_j$'s, i.e. $
m^{\underline\sigma}_N=1/N\,\sum_{j=1}^N \sigma_j$,
the mean choice of the perceived choices of all individuals.
\end{enumerate}
In \cite{DPRST09} we consider a system comprised by $N$ identical 
individuals. The evolution is characterized by two parameters $\b$ and 
$\g$, where $\b$ determines the randomness in the perception process, 
while $\g$ determines the randomness in the process of adaptation to the 
mean of other individuals. It is a parameter that expresses a measure of the disutility of non-conformance. When $\gamma >0$ individuals tend to conform their behavior to the mean behavior of all the others; conversely, when $\gamma <0$, there is an incentive to non-conformity. We take into account only positive $\gamma$, since, as shown in \cite{BlDu03} and in \cite{BrDu01}, the results are more interesting both from the social and technical point of view than in the other case.\\ 
For simplicity, we consider a binary decision problem for individual agents, so both $\o_i$ and $\s_i$ are 
assumed to take values $\pm 1$, and we call them {\em spins}. The main 
results in \cite{DPRST09} concern the limiting dynamics as $N \ra 
+\infty$, and normal fluctuations about this limit. In particular we 
show that the parameter space is divided into two main regions, the {\em subcritical} and the {\em supercritical} ones, corresponding respectively to a minor or a major incentive to conform. For all value of the parameters, a ``neutral'' equilibrium solution  exists, and, roughly speaking, it corresponds to equal proportion of individuals with state $+1$ and $-1$.  In the supercritical region,  however, other equilibria appear, with a strict majority of  individuals in one of the two states, as effect of a stronger attitude  to conformism.

%The aim of this paper is twofold. On one hand we weaken the homogeneity 
%assumption  of identical individuals; on the other we obtain a scaling 
%limit of the dynamics of fluctuations for {\em critical} values of the 
%parameters, i.e. in the boundary between the subcritical and the 
%supercritical region.

In this paper we weaken the homogeneity assumption of identical individuals; here individuals are divided into reference groups. The belonging to a given group is coded in a parameter, that we assume {\em random} but {\em constant in time}. By adopting the terminology used in Statistical Mechanics of {\em disordered systems}, we will refer to the set of these parameters as {\em random field}. Depending on the particular application, the random field may describe affiliation of individuals to different social classes, ethnic or religious groups, geographic areas and so on. In a financial setting, one can model firms of different dimension, or acting in different markets. \\
Our aim is twofold.
On one hand, following the same approach used in \cite{DPRST09} to study the limiting dynamics as $N \ra +\infty$, we show that the phase diagram of the model is more complex than the one of the homogeneous case. \\
On the other we obtain a scaling limit of the dynamics of fluctuations for {\em critical} values of the  parameters, i.e. in the boundary between the subcritical and the supercritical region.

\noindent
%We divide the individuals into groups, each group evolving with 
%different dynamic rules, by introducing a random field which assigns to 
%each individual the group he belongs to. 
We consider the case of two 
groups, i.e. a $\pm 1$-valued random field. This greatly simplifies the 
analysis of critical fluctuations, since it allows the reduction to a 
low dimensional order parameter. General random fields would require the 
spectral analysis of operators that, by the non reversibility of the 
system, are not self-adjoint. This difficulty is innocuous when 
reduction to low dimension is possible.

\noindent
%In order to illustrate critical fluctuations, 
To illustrate our main results, we need to introduce some 
notations. For times $t \geq 0$, let $\s_i(t),\o_i(t) \in \{-1,+1\}$ be 
the values of the spins at time $t$, which evolve as a continuous-time 
Markov chain. Moreover, let $\eta_i \in \{-1,+1\}$ be the group of 
the $i$-th individual. Denote by
\[
\rho_N(t) := \frac{1}{N} \sum_{i=1}^N \delta_{(\s_i(t),\o_i(t),\eta_i)}
\]
the empirical measure at time $t$. In the model we introduce in Section 
2, $\rho_N(t)$ evolves itself as a Markov process (this fact could be 
indeed taken as definition of mean field dynamics). $\rho_N$ is a 
probability on $\{-1,+1\}^3$, so it lives on a linear manifold of 
dimension $7$. As we shall see, many explicit computations are made 
possible by this low dimensionality, that would be lost in more general 
cases, for instance when $\eta_i$ could take infinitely many values. As 
$N \ra +\infty$, $\rho_N(t)$ converges in probability to a deterministic 
flow $\rho(t)$, which is a solution of an ordinary differential equation 
(ODE). This result corresponds to a \emph{law of large numbers} and it allows us to describe the macroscopic evolution of the system, which is deterministic. We provide the full phase diagram of the stationary solutions of 
this ODE, in terms of the parameters of the model. \\
Since in real systems $N$ is large but finite, it is relevant to obtain first order (Normal) corrections to the $N \ra +\infty$ limiting dynamics.
So, we show that the 
{\em fluctuation process}
\[
\sqrt{N} [ \rho_N(t) - \rho(t) ]
\]
converges in law to a Gaussian process (\emph{Central Limit Theorem}), whose covariance is determined explicitely. Gaussian approximations of empirical mean are widely used in applications. For example, in the financial context of {\em credit risk analysis} (see e.g. \cite{GiWe04,GiWe06, DPRST09,DaPTo}),  states are indicators of the financial health of a firm; Gaussian approximations allow to compute quantiles of the excess losses suffered by a financial institution holding a large portfolio with positions issued by the firms. Moreover, whenever parameters have to be estimated from data, the Central Limit Theorem yields asymptotic normality of estimators, which is in practice a very desirable property.
\\
In the subcritical region, where only a ``neutral'' equilibrium exists, it can be shown that the estimates of Normal fluctuations are uniform in time. When the parameters approach the critical values, these estimates loose their accurancy and we need a different scaling to describe better their behavior. In fact, for critical values of the 
parameters, it is expected that long time fluctuations are such that a 
space-time scaling of the form
\begin{equation}
\label{intro1}
N^{1/4} [\rho_N(N^{\alpha} t) - \overline{\rho}],
\end{equation}
where $\overline{\rho} := \lim_{t \ra +\infty} \rho(t)$, has a 
nontrivial limit in law. Although many homogeneous and reversible models 
are well understood in this respect, this paper provides, to our 
knowledge, the first example of study of critical fluctuations for a non 
reversible, inhomogeneous model. We show that, in the homogeneous case 
($\eta_i \equiv 0$), the ``standard'' $\alpha = 1/2$ scaling holds true, 
and (\ref{intro1}) converges to a cubic diffusion. When instead the 
$\eta_i$'s are i.i.d. and nonzero, the space fluctuations of the field 
destroy the above picture: critical fluctuations appear at a much 
shorter time $\alpha = 1/4$, and are driven by the normal fluctuations 
of the field. In some sense, this is a dynamic analog of the result in 
\cite{MaPe91} for the fluctuation of the Curie-Weiss model at critical 
temperature. We remark that the dynamics of critical fluctuations are known to exhibit universality features (see \cite{Da83}). Our result shows that the presence of inhomogeneities may lead to a new universality class (see \cite{CoDa10} for more results on this subject).

The paper is organized as follows. In Section 2 we define the class of 
models we study; Section 3 contains the results concerning the limiting 
dynamics, phase diagram and normal fluctuations; Section 4 is devoted to 
critical fluctuations; proofs of all results are finally given in Section 5.

%\subsection{General aspects}

%\subsection{Purpose and modeling aspects}

%\subsection{Financial application}

%\subsection{Methodology}\label{meth}

\section{The model}

\subsection{Description of the model}

Let $\mathscr{S}=\{-1,+1\}$ and $\underline{\eta}=(\eta_j)_{j=1}^N \in \mathscr{S}^N$ be a sequence of independent, identically distributed, symmetric, Bernoulli random variables defined on some probability space $(\Omega, \mathcal{F}, P)$, that is $P(\eta_j=-1)=P(\eta_j=+1)=1/2$, for any $j$. We indicate by $\mu$ their common law.\\
Given a configuration $(\underline{\sigma}, \underline{\omega})=(\sigma_j, \omega_j)_{j=1}^N \in \mathscr{S}^{2N} $ and a realization of the random medium $\underline{\eta}$, we construct a $2N$-spin system evolving as a continuous-time Markov chain on $\mathscr{S}^{2N}$, with infinitesimal generator $L_N$ acting on functions 
$f:\mathscr{S}^{2N} \longrightarrow \mathbb{R}$ as follows:
\begin{equation}\label{IG1'}
    L_Nf(\underline{\sigma},\underline{\omega})=\sum_{j=1}^{N} e^{-\beta \sigma_j \omega_j} \nabla_j^\sigma f(\underline{\sigma}, \underline{\omega}) + \sum_{j=1}^{N} e^{-\gamma \omega_j( m^{\underline{\sigma}}_{N}+ h \eta_j)} \nabla_j^\omega f(\underline{\sigma}, \underline{\omega}),
\end{equation}
where $\nabla_j^\sigma f(\underline{\sigma}, \underline{\omega})=f(\underline{\sigma}^j, \underline{\omega})-f(\underline{\sigma}, \underline{\omega})$ and $\nabla_j^\omega f(\underline{\sigma}, \underline{\omega})=f(\underline{\sigma}, \underline{\omega}^j)-f(\underline{\sigma}, \underline{\omega})$. The $k$-th component of $\underline{\sigma}^j$, which has the meaning of a $\sigma$-spin flip at site $j$, is
\begin{displaymath}
    \sigma^j_k=\left\{
    \begin{array}{rcc}
        \sigma_k & \mathrm{for} &k\neq j\\
        -\sigma_k& \mathrm{for} &k= j
    \end{array}
    \right.
\end{displaymath}
and the $\omega$-spin flip at site $j$ is defined similarly. The parameters $\beta$, $\gamma$ and $h$ are positive. \\
The quantities \mbox{$e^{-\beta\sigma_j\omega_j}$} and \mbox{$e^{-\gamma \omega_j( m^{\underline{\sigma}}_{N}+ h \eta_j)}$} represent the jump rates of the spins, the rates at which the transitions $\sigma_j\longrightarrow -\sigma_j$ and  $\omega_j\longrightarrow -\omega_j$  occur respectively for some $j$. The rates $e^{- \beta \sigma_j \omega_j}$ describe how states are perceived: for $\beta = 0$ perceived states are completely random, while alignment to the real states improves as $\beta$ grows. The rates $e^{-\gamma \omega_j( m^{\underline{\sigma}}_{N}+ h \eta_j)}$ are comprised by two factors: $e^{-\gamma \omega_j m^{\underline{\sigma}}_{N}}$ incentivize alignment with the average perceived state of the community ({\em conformism}), while the factor $e^{-\gamma h \omega_j \eta_j}$ models different attitudes within different reference groups.

The expression \eqref{IG1'} describes a system of mean field coupled pairs of spins, each with its own random environment. It is subject to an inhomogeneous interaction (of intensity $h$) parametrized by the components $\eta_j$. With the expression mean field we mean that the sites interact all each other in the same way and this assumption allows us to suppose that the interaction depends on the value of the magnetization
\[
    m^{\underline{\sigma}}_{N}(t)=\frac{1}{N}\sum_{j=1}^N\sigma_j(t).
\]
The initial condition $(\underline{\sigma}(0),\underline{\omega}(0))$ is assumed to have product distribution $\lambda^{\otimes N}$, where $\lambda$ is a probability measure on $\mathscr{S}^2$.\\
The quantity $(\sigma_j(t),\omega_j(t))$ represents the time evolution on $[0,T]$, $T$ fixed, of $j$-th pair of spin values; it is the trajectory of the single $j$-th pair of spin values in time. The space of all these paths is $(\mathcal{D}[0,T])^2$, where $\mathcal{D}[0,T]$ is the space of the right-continuous, piecewise-constant functions from $[0,T]$ to $\mathscr{S}$, endowed with the Skorohod topology, which provides a metric and a Borel $\sigma$-field (as we can see in \cite{EtKu86}).

\section{Approach and main results}

\subsection{Deterministic limit: law of large numbers}

The operator $L_N$ given in \eqref{IG1'} defines an irreducible, finite-state Markov chain. It follows that the process admits a unique stationary distribution $\nu_N$, but it can be proved that our model is non reversible (see the analogous proof for the homogeneous model in \cite{DPRST09}). This fact implies that an explicit formula for the stationary distribution $\nu_N$ and its $N \longrightarrow +\infty$ asymptotics are not available. So, we follow a dynamic approach. This means that first, we derive the dynamics of the process \eqref{IG1'}, in the limit as $N\longrightarrow +\infty$, in a fixed time interval $[0,T]$ and later, we study the large time behavior of the limiting dynamics. %This is not necessarily equivalent to studying the $N \longrightarrow +\infty$ properties of $\nu_N$.  However, as we shall show later in this paper, this provides rather sharp information on how the system behaves for $t$ and $N$ large.

So, let $(\sigma_j [0,T],\omega_j[0,T])_{j=1}^N \in (\mathcal{D}[0,T])^{2N}$ denote a path of the system in the time interval $[0,T]$, with $T$ positive and fixed. If $f(\sigma_j [0,T], \omega_j[0,T])$ is a function of the trajectory of a single pair of spins, we are interested in the asymptotic behavior of \emph{empirical averages} of the form
\[ 
\frac{1}{N} \sum_{j=1}^{N} f(\sigma_j [0,T], \omega_j [0,T]) =: \int f d\rho_N \,,
\]
where $\{\rho_N\}_{N \geq 1}$ is the sequence of \emph{empirical measures}
\[ 
\rho_N := \frac{1}{N} \sum_{j=1}^{N} \delta_{(\sigma_j [0,T], \omega_j [0,T], \eta_j)}\,. 
\]
We may think of $\rho_N$ as a random element of $\mathcal{M}_1((\mathcal{D}[0,T])^2 \times \mathscr{S})$, the space of probability measures on $(\mathcal{D}[0,T])^2 \times \mathscr{S}$ endowed with the weak convergence topology. 

First, we want to determine the weak limit of $\rho_N$ in $\mathcal{M}_1((\mathcal{D}[0,T])^2 \times \mathscr{S})$, when $N$ grows to infinity, i.e. for $f \in \mathcal{C}_{b}$ we look for $\lim_{N \rightarrow +\infty} \int f d\rho_N$. It corresponds to a law of large numbers, where the limit is a deterministic measure. Being an element of $\mathcal{M}_1((\mathcal{D}[0,T])^2 \times \mathscr{S})$, such a limit can be viewed as a stochastic process, which represents the dynamics of the system in the infinite volume limit.

The result we actually present is a \emph{large deviation principle}, which is much stronger than a law of large numbers. We start with some preliminary notions letting, in what follows, $W \in \mathcal{M}_1((\mathcal{D}[0,T])^2)$ denote the law of the $\mathscr{S}^2$-valued process $(\sigma(t), \omega(t))_{t \in [0,T]}$, such that the initial condition $(\sigma(0), \omega(0))$ has distribution $\lambda$ and both $\sigma(\cdot)$ and $\omega(\cdot)$ change sign with constant rate equal to 1. 
By $W^{\otimes N}$ we mean the product of $N$ copies of $W$, which represents the law of the $2N$-spin system %, whose generator is \eqref{IG1'}, where we have set $c^{\underline{\eta}, \, \underline{\sigma}}_N = c^{\underline{\eta}, \, \underline{\omega}}_N \equiv 1$, in other words, the law of our system
 in absence of interaction. Moreover, we shall denote by $P_N^{\underline{\eta}}$ the law of the process $(\underline{\sigma}([0,T]), \underline{\omega}([0,T])) = (\underline{\sigma}(t), \underline{\omega}(t))_{t \in [0,T]}$, with infinitesimal generator \eqref{IG1'} and initial distribution $\lambda^{\otimes N}$, for a given $\underline{\eta}$.

For $Q \in \mathcal{M}_1((\mathcal{D}[0,T])^2 \times \mathscr{S})$, let 
\[
H(Q \vert W \otimes \mu) := \left \{
\begin{array}{ll}
\int dQ \log \frac{dQ}{d(W \otimes \mu)} & \mbox{if} \quad Q \ll W \otimes \mu \quad \mbox{and} \quad \log \frac{dQ}{d(W \otimes \mu)} \in L^1(Q)\,,\\
+\infty & \mbox{otherwise}\,,
\end{array}
\right.
\]
denote the \emph{relative entropy} between $Q$ and $W\otimes\mu$. Moreover, $\Pi_t Q$ denotes the marginal law of $Q$ at time $t$, and
\[
m^{\sigma}_{\Pi_tQ} := \int_{\mathscr{S}^3} \sigma \Pi_t Q(d\sigma, d\omega, d\eta)\,.
\]
For a given path $(\sigma([0,T]), \omega([0,T])) \in (\mathcal{D}[0,T])^2$, let $\mathcal{N}_t^{\sigma}$ (resp. $\mathcal{N}_t^{\omega}$) be the process counting the jumps of $\sigma(\cdot)$ (resp. $\omega(\cdot)$). Define 
\begin{allowdisplaybreaks}
\begin{multline}\label{F(Q)'}
F(Q)  := \int  \bigg[  \int_0^T  \left( 1 - e^{-\beta \sigma(t) \omega(t)}\right)dt + \beta  \int_0^T  \sigma(t) \omega(t) d\mathcal{N}^{\sigma}_t\\
         \quad + \int_0^T  \left( 1 - e^{-\gamma \omega(t) \left( m^{\sigma}_{\Pi_tQ} +  h \eta \right)}\right)dt + \gamma  \int_0^T  \omega(t) \left( m^{\sigma}_{\Pi_tQ} +  h \eta\right) d\mathcal{N}^{\sigma}_t  \bigg]dQ \,,
\end{multline}
\end{allowdisplaybreaks}\\
whenever 
$$\int \left(\mathcal{N}_T^{\sigma} + \mathcal{N}_T^{\omega} \right) dQ < +\infty\,,$$
and $F(Q)=0$ otherwise. Finally, let
\[
I(Q) := H(Q \vert W \otimes \mu) - F(Q)\,.
\]
We remark that, if $\int (\mathcal{N}_T^{\sigma} + \mathcal{N}_T^{\omega})dQ=+\infty$, then $H(Q \vert W \otimes \mu)=+\infty$ (see Lemma 5.4 in \cite{DPRST09}) and thus also $I(Q)=+\infty$.

\begin{PROP}\label{LDP}
%For each $Q \in {\cal{M}}_1 (\mathcal{D}[0,T])^2$, $I(Q) \geq 0$, and $I(\cdot)$ is a
%lower-semicontinuous function with compact level-sets (i.e. for each
%$k>0$ one has that $\{Q:I(Q) \leq k\}$ is compact in the weak
%topology). Moreover, for $A,C \subseteq {\cal{M}}_1
%(\mathcal{D}[0,T])^2$ respectively open and
%closed for the weak topology, we have
%\begin{eqnarray}
%\liminf_N \frac{1}{N} \log P(\rho_N \in A) & \geq & - \inf_{Q \in A} I(Q) \label{lb}\\
%\limsup_N \frac{1}{N} \log P(\rho_N \in C) & \leq & - \inf_{Q \in C} I(Q). \label{ub}
%\end{eqnarray}
%This means that 
The distributions of $\rho_N$ obey a {\em large
deviation principle} (LDP) with {\em rate function} $I(\cdot)$ (see
e.g. \cite{DeZe93} for the definition and fundamental facts on LDP).
\end{PROP}

The key step to derive a law of large numbers from Proposition \ref{LDP} is given in the following result.

\begin{PROP}\label{thmMKV'}
The equation $I(Q)=0$ has a unique solution $Q_* \in \mathcal{M}_1 ((\mathcal{D}[0,T])^2 \times \mathscr{S})$ which admits the decomposition $Q_*(d\s[0,T], d\omega[0,T], d\eta) = Q_*^{\eta}(d\s[0,T], d\omega[0,T]) \mu(d\eta)$. Moreover, the marginals \mbox{$q_t^{\eta}=\Pi_t Q_*^{\eta} \in \mathcal{M}_1 (\mathscr{S}^2)$} are weak solutions of the nonlinear McKean-Vlasov equation
\begin{equation}\label{MKV1'}
    \left\{
    \begin{array}{cccr}
        \frac{\partial q_t^{\eta}}{\partial t} & = & \mathcal{L}^{\eta} q_t^{\eta} & \qquad (t \in [0,T], \; \eta \in \mathscr{S})\\
        q_0^{\eta} & = & \lambda & \\
    \end{array}
    \right.\,,
\end{equation}
where, for all the triples  $(\sigma, \omega, \eta) \in \mathscr{S}^3$, the operator $\mathcal{L}^{\eta}$ acts as follows:
\begin{equation}\label{IG2'}
    \mathcal{L}^{\eta} q_t^{\eta} (\sigma, \omega)=\nabla^\sigma \left[ e^{-\beta \sigma \omega} q_t^{\eta} (\sigma, \omega) \right] + \nabla^\omega \left[ e^{-\gamma \omega ( m^{\sigma}_{q_t} +  h \eta)} q_t^{\eta} (\sigma, \omega) \right]\,,
\end{equation}
and $q_t$ is defined by 
\[q_t(\sigma, \omega) = \int_{\mathscr{S}} q_t^{\eta} (\sigma, \omega) \mu(d\eta)\,.\]
\end{PROP}

From Propositions \ref{LDP} and \ref{thmMKV'}, it is easy to derive the following \emph{strong law of large numbers}.

\begin{TEO}\label{lln}
Let $Q_* \in \mathcal{M}_1 ((\mathcal{D}[0,T])^2 \times \mathscr{S})$ be the probability given in Proposition \ref{thmMKV'}. Then 
$$\rho_N\longrightarrow Q_*\qquad\text{almost surely}$$ 
in the weak topology.
\end{TEO}

\begin{REM}
\label{propchaos}
The result in Theorem \ref{lln} shows the convergence of the sequence of the empirical measures. The qualitative and quantitative analysis of its limit will be treated in the next section. It is worth to point out here a consequence of Theorem \ref{lln}. Let $i_1, i_2, \ldots i_m$ be fixed indexes in $\{1,2,\ldots,N\}$. Then,  the joint law of the random variables $(\s_{i_j}[0,T], \o_{i_j}[0,T], \eta_{i_j})_{j=1}^m$ converges weakly to $Q_*^{\otimes m} $. This can be shown along the same lines of Theorem 3 in \cite{CoEi88}, and it is known as {\em propagation of chaos} property: in the limit as $N \ra +\infty$, the joint law of the state evolutions of given  individuals is a product measure; moreover, single individuals evolve their state according to the law $Q_*$.
\end{REM}

The proofs of Propositions \ref{LDP} and \ref{thmMKV'} and of Theorem \ref{lln} are based on large deviations techniques applied to mean field models, first introduced in \cite{DaPdHo95} and then generalized in \cite{DPRST09} for non reversible mean field models. They present various technical difficulties due to the
unboundedness and non continuity of $F$, which are related to the non reversibility of
the model. They are not given here, because they are analogous to the same results in the homogeneous case (see Proposition 3.1, Proposition 3.2 and Theorem 3.3 in \cite{DPRST09}).  

\subsection{Equilibria of the limiting dynamics: phase
transition}
%stationary solutions

The equation \eqref{MKV1'} describes the behavior of the system governed by generator \eqref{IG1'} in the infinite volume limit, i.e. of infinitely many individuals. We are interested in the detection of the \hbox{$t$-stationary} solution(s) of this equation and in the study of its (their) stability properties. We recall that, to be $t$-stationary solution(s) for \eqref{MKV1'} it has to be satisfied the equation $\mathcal{L}^{\eta} q^{\eta}=0$, for every $t$.

First of all, we proceed to reformulate the ``original'' McKean-Vlasov equation \eqref{MKV1'} in terms of $m^{\eta}_{q_t}$, $m^{\sigma}_{q_t}$, $m^{\omega}_{q_t}$, $m^{\sigma \omega}_{q_t}$, $m^{\sigma \eta}_{q_t}$, $m^{\omega \eta}_{q_t}$, $m^{\sigma \omega \eta}_{q_t}$ defined as follows:
\begin{equation} \label{mag2d'}
m^{\eta}_{t} := \frac{1}{2} \sum_{\sigma, \omega \in \mathscr{S}} \sum_{\eta \in \mathscr{S}} \eta \, q_t^{\eta}(\sigma, \omega)  
\end{equation}
\begin{allowdisplaybreaks}
\begin{align}
\label{mag2a'} m^{\sigma}_{t} &:= \frac{1}{2} \sum_{\sigma, \omega \in \mathscr{S}} \sum_{\eta \in \mathscr{S}} \sigma \, q_t^{\eta}(\sigma, \omega)  \qquad m^{\sigma \eta}_{t} := \frac{1}{2} \sum_{\sigma, \omega \in \mathscr{S}} \sum_{\eta \in \mathscr{S}} \sigma \eta \, q_t^{\eta}(\sigma, \omega) \\
\label{mag2b'} m^{\omega}_{t} &:= \frac{1}{2} \sum_{\sigma, \omega \in \mathscr{S}} \sum_{\eta \in \mathscr{S}} \omega \, q_t^{\eta}(\sigma, \omega) \qquad m^{\omega \eta}_{t} := \frac{1}{2} \sum_{\sigma, \omega \in \mathscr{S}} \sum_{\eta \in \mathscr{S}} \omega \eta \, q_t^{\eta}(\sigma, \omega)  \\
\label{mag2c'} m^{\sigma \omega}_{t} &:= \frac{1}{2} \sum_{\sigma, \omega \in \mathscr{S}} \sum_{\eta \in \mathscr{S}} \sigma \omega \, q_t^{\eta}(\sigma, \omega)\qquad m^{\sigma \omega \eta}_{t} := \frac{1}{2} \sum_{\sigma, \omega \in \mathscr{S}} \sum_{\eta \in \mathscr{S}} \sigma \omega \eta \, q_t^{\eta}(\sigma, \omega) \,,
\end{align}
\end{allowdisplaybreaks}\\
where $q_t^{\eta}$ has the meaning explained in Proposition \ref{thmMKV'} and we have written $m^{\cdot}_{t}$ instead of $m^{\cdot}_{q_t}$. We introduce these expectations, because the probability measure $q_t$ on $\mathscr{S}^3$ is completely determined by them.\\
The quantities defined above, or simple functions of them, have natural interpretations in the social setting. For a given time $t$, $m_t^\sigma$ (resp. $m_t^\omega$) clearly represents the averaged perceived (resp. real) state. Then, $\frac{1+m^{\sigma \omega}_t}{2}$ is the probability that the perceived state ($\sigma$) of an individual is equal to its real state ($\omega$). Also, $m^{\omega}_t \pm m^{\omega \eta}_t$ is the average state {\em within} the reference group associated to $\eta = \pm 1$. Similar interpretations can be given for the other quantities in (\ref{mag2d'})-(\ref{mag2c'}).

\begin{LEMMA}\label{MKV2'}
Equations \eqref{MKV1'} can be rewritten in the following form:
%\begin{allowdisplaybreaks}
%\begin{align*}
\begin{equation}\label{diffeq}
\begin{split}
&\dot{m}^{\eta}_{t} \;\;\,=\, 0 \\
&\dot{m}^{\sigma}_{t}\;\;\, = -2\, m^{\sigma}_{t} \cosh(\beta ) + 2 \, m^{\omega}_{t} \sinh(\beta ) \\
&\dot{m}^{\omega}_{t} \;\;\,= -2 \, m^{\omega}_{t} \cosh(\gamma h) \cosh(\gamma m^{\sigma}_{t}) - 2 \, m^{\omega \eta}_{t} \sinh(\gamma h) \sinh(\gamma m^{\sigma}_{t})\\
                                &\qquad\qquad + 2 \, \cosh(\gamma h) \sinh(\gamma m^{\sigma}_{t})\\ 
&\dot{m}^{\sigma \omega}_{t} \;=\, 2 \, m^{\sigma}_{t} \cosh(\gamma h) \sinh(\gamma m^{\sigma}_{t}) - 2 \, m^{\sigma \omega}_{t} \left[ \cosh(\beta) + \cosh(\gamma h) \cosh(\gamma m^{\sigma}_{t}) \right]\\
                                            &  \qquad\qquad + 2 \, m^{\sigma \eta}_{t} \sinh(\gamma h) \cosh(\gamma m^{\sigma}_{t}) - 2 \, m^{\sigma \omega \eta}_{t} \sinh(\gamma h) \sinh(\gamma m^{\sigma}_{t})+2 \sinh(\beta) \\ 
&\dot{m}^{\sigma \eta}_{t} \;\,=-2\, m^{\sigma \eta}_{t} \cosh(\beta) + 2\, m^{\omega \eta}_{t} \sinh(\beta) \\
&\dot{m}^{\omega \eta}_{t} \;\,= -2 \,  m^{\omega}_{t} \sinh(\gamma h) \sinh(\gamma m^{\sigma}_{t}) - 2 \, m^{\omega \eta}_{t} \cosh(\gamma h) \cosh(\gamma m^{\sigma}_{t})\\
                                       &  \qquad\qquad + 2 \, \sinh(\gamma h) \cosh(\gamma m^{\sigma}_{t}) \\
&\dot{m}^{\sigma \omega \eta}_{t} =\,2 \, m^{\sigma}_{t} \sinh(\gamma h) \cosh(\gamma m^{\sigma}_{t}) - 2 \, m^{\sigma \omega}_{t} \sinh(\gamma h) \sinh(\gamma m^{\sigma}_{t})\\
                                      &  \qquad\qquad + 2 \, m^{\sigma \eta}_{t} \cosh(\gamma h) \sinh(\gamma m^{\sigma}_{t}) - 2 \, m^{\sigma \omega \eta}_{t} \left[ \cosh(\beta) + \cosh(\gamma h) \cosh(\gamma m^{\sigma}_{t}) \right]\,,
 \end{split}
 \end{equation}                                 
%\end{align*}
%\end{allowdisplaybreaks}\\
with initial condition $m^{\eta}_{0} = m^{\eta}_{\lambda} = 0$, $m^{\sigma}_{0}= m^{\sigma}_{\lambda}$, $m^{\omega}_{0}= m^{\omega}_{\lambda}$, $m^{\sigma \omega}_{0}= m^{\sigma \omega}_{\lambda}$, $m^{\sigma \eta}_{0}= m^{\sigma \eta}_{\lambda}$, $m^{\omega \eta}_{0}= m^{\omega \eta}_{\lambda}$ and $m^{\sigma \omega \eta}_{0}= m^{\sigma \omega \eta}_{\lambda}$.
\end{LEMMA}

\Proof See Section \ref{proofs}.
\fine

The variable $m^{\eta}_t$ is static, thus any equilibrium solution of the system in Lemma \ref{MKV2'} is of the form
\begin{equation}\label{soleq'}
\begin{split}
 & m^{\sigma}_{*}\;\;\, = \tanh(\beta) \frac{\sinh(\gamma m^{\sigma}_{*}) \cosh(\gamma m^{\sigma}_{*})}{\cosh^2(\gamma  m^{\sigma}_{*}) + \sinh^2(\gamma h)}\\
& m^{\omega}_{*} \;\;\,= \frac{\sinh(\gamma m^{\sigma}_{*}) \cosh(\gamma m^{\sigma}_{*})}{\cosh^2(\gamma m^{\sigma}_{*}) + \sinh^2(\gamma h)}\\
& m^{\sigma \omega}_{*}\; = \dots\\
& m^{\sigma \eta}_{*} \;\,= \tanh(\beta) \tanh(\gamma h) \frac{1 + \sinh^2(\gamma h)}{\cosh^2(\gamma m^{\sigma}_{*}) + \sinh^2(\gamma h)}\\
& m^{\omega \eta}_{*} \;\,=\tanh(\gamma h) \frac{1 + \sinh^2(\gamma h)}{\cosh^2(\gamma m^{\sigma}_{*}) + \sinh^2(\gamma h)}\\
& m^{\sigma \omega \eta}_{*} = \dots\,. 
\end{split}
\end{equation}

To discover the presence of phase transition(s) (multiple equilibria) and the stability of equilibria, it is sufficient studying the first equation of \eqref{soleq'}:
\[ 
m^{\sigma}_{*} = \tanh(\beta) \frac{\sinh(\gamma m^{\sigma}_{*}) \cosh(\gamma m^{\sigma}_{*})}{\cosh^2(\gamma m^{\sigma}_{*}) + \sinh^2(\gamma h)}\,,
\]
because all the remaining are \mbox{$ m^{\cdot}_{*}= m^{\cdot}_{*}( m^{\sigma}_{*})$}, hence $\displaystyle{\lim_{t\rightarrow+\infty} m^{\cdot}_{t}= m^{\cdot}_{*}}$, when $\displaystyle{\lim_{t\rightarrow+\infty} m^{\sigma}_{t}= m^{\sigma}_{*}}$. The stationary system we are dealing with is essentially one-dimensional.
 
\begin{figure}[htb]
\centering%
\subfigure[]{\includegraphics[height=5cm,width=0.49\textwidth]{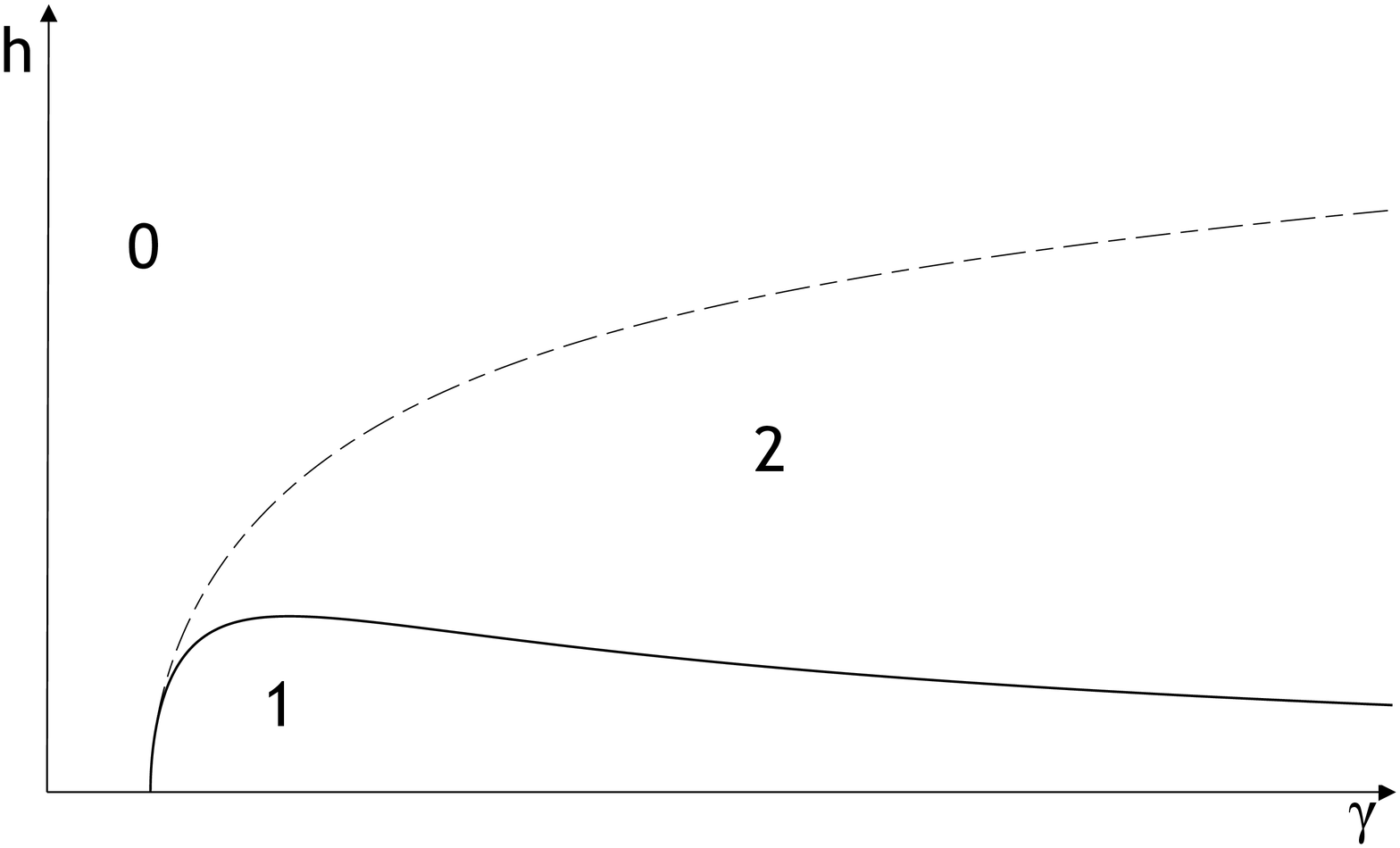}\label{fig2a}} 
\subfigure[]{\includegraphics[height=5cm, width=0.49\textwidth]{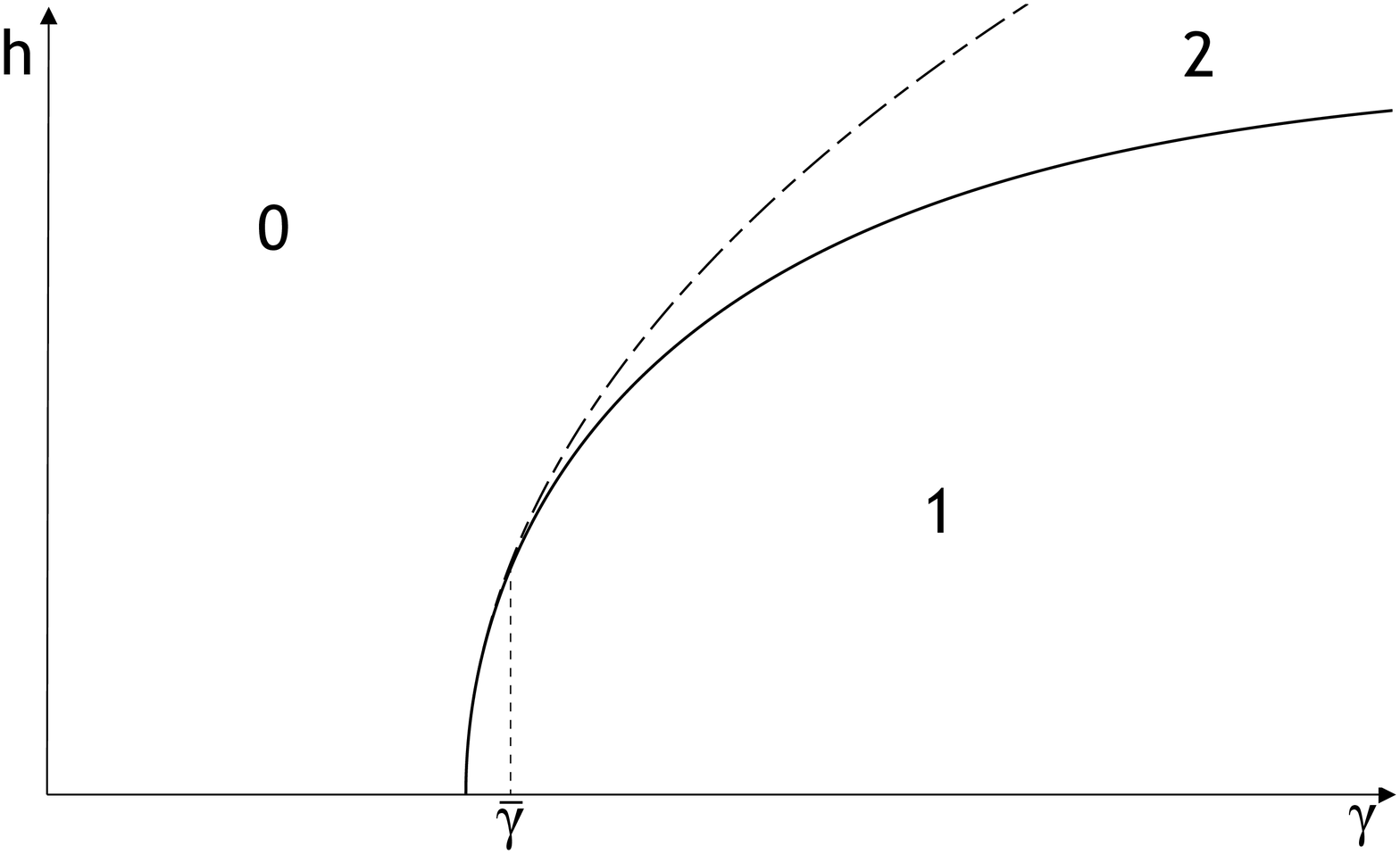}\label{fig2b}}
%\subfigure[]{\includegraphics[height=5cm,width=0.49\textwidth]{PD.pdf}\label{fig2a}} 
%\subfigure[]{\includegraphics[height=5cm, width=0.49\textwidth]{PDz.pdf}\label{fig2b}}
\caption{Phase diagram for a fixed value of $\beta$ (figure \ref{fig2a})  and zoom of the area where the bifurcation occurs (figure \ref{fig2b}).}
\label{fig:2}
\end{figure}
%Se volete togliere le didascalie "personalizzate" e la sottonumerazione per le due figure bisogna togliere le parentesi quadre [caption]. Se volete avere la sottonumerazione, ma non le didascalie, basta lasciare vuote tali parentesi.

For a fixed value of $\beta$, the phase diagram is qualitatively drawn in Figure \ref{fig:2}. There are three phases, corresponding to 0, 1 and 2  solutions having $m^{\sigma}_* >0$, respectively. By symmetry, we have the same number of solutions with $m^{\sigma}_* <0$.\\ 
The continuous separation curve is
\begin{equation}\label{curve'}
h = h(\beta, \gamma) = \frac{1}{\gamma} \mathrm{arccosh} (\sqrt{\gamma \tanh(\beta)})\,, \qquad \gamma \in \left[\frac{1}{\tanh(\beta)}, +\infty \right) \,,
\end{equation}
while the dotted one is obtained numerically and it is due to the fact that the function
\begin{equation}\label{Gamma}
m^{\sigma}_{*} \longmapsto \Gamma_{\beta, \gamma, h} (m^{\sigma}_{*}) := \tanh(\beta) \frac{\sinh(\gamma m^{\sigma}_{*}) \cosh(\gamma m^{\sigma}_{*})}{\cosh^2(\gamma m^{\sigma}_{*}) + \sinh^2(\gamma h)}
\end{equation}
is not always concave.\\ The two curves coincide for $\gamma \in \left[\frac{1}{\tanh(\beta)}, \frac{3}{2\tanh(\beta)}\right]$ and separate at the ``tricritical'' point $(\bar{\gamma}, h\left(\beta,\bar{\gamma}\right))=\left(\frac{3}{2\tanh(\beta)}, h \left(\beta, \frac{3}{2\tanh(\beta)}\right)\right)$.

\begin{TEO}\label{solutions}
Consider the equations \eqref{soleq'} and fix a value for $\beta$. The point
\[
\underline m_{*}^{0}  :=  \left( 0, 0, \frac{\mathrm{th}(\beta) \mathrm{th}(\gamma h) \mathrm{sh}(\gamma h) + \mathrm{sh}(\beta)}{\mathrm{ch}(\beta) + \mathrm{ch}(\gamma h)}, \mathrm{th}(\beta) \mathrm{th}(\gamma h), \mathrm{th}(\gamma h), 0  \right) \,.
\] 
is a solution for all values of the parameters. 
\begin{enumerate}
\item
If $\gamma < \frac{\cosh^2(\gamma h)}{\tanh(\beta)}$ (Phases \emph{\textsf{0}}  and \emph{\textsf{2}}  in  Figure \ref{fig:2}), then $\underline m_{*}^{0}$ is linearly stable for (\ref{diffeq}).
\item
If $\gamma < \frac{\cosh^2(\gamma h)}{\tanh(\beta)}$, and $h$ is above the dotted line in Figure \ref{fig:2}, which has been obtained numerically (Phase \emph{\textsf{0}}),  then $\underline m_{*}^{0}$ is the unique solution of \eqref{soleq'}.
\item
If $\gamma > \frac{\cosh^2(\gamma h)}{\tanh(\beta)}$, i.e. $(\gamma,h)$ is below the curve \eqref{curve'} , the continuous curve  in Figure \ref{fig:2} (Phase \emph{\textsf{1}}), then  \eqref{soleq'} has three solutions: \\$\underline m_{*}^{0}$,\hfill $(m_{*}, m^{\omega}_{*}(m_{*}), m^{\sigma \omega}_{*}(m_{*}), m^{\sigma \eta}_{*}(m_{*}), m^{\omega \eta}_{*}(m_{*}), m^{\sigma \omega \eta}_{*}(m_{*}))$\hfill and\hfill
$(-m_{*}, -m^{\omega}_{*}(m_{*}),$\par $m^{\sigma \omega}_{*}(-m_{*}), m^{\sigma \eta}_{*}(m_{*}), m^{\omega \eta}_{*}(m_{*}), m^{\sigma \omega \eta}_{*}(-m_{*}))$, where $m_{*}$ is the unique positive solution of the first equation of \eqref{soleq'}. 

\item If we choose the parameters above the curve \eqref{curve'} and $h$ is small enough, in other words if $(\gamma, h)$ belongs to the Phase \emph{\textsf{2}} of Figure \ref{fig:2}, then two further solutions arise.
\end{enumerate}
\end{TEO}

\Proof See Section \ref{proofs}.
\fine

Interpretations of Theorem \ref{solutions} and of the phase diagram in Figure \ref{fig:2} are most easily given in terms of {\em opinion dynamics}, where the states of individuals are (binary) opinions on a given subject. Phase \textsf{0} in Figure \ref{fig:2} can be seen as a small perturbation of the case where individuals choose their opinion randomly, and independently of the others. \\
For small inhomogeneity ($h$ small), as $\gamma$ crosses the curve $h = h(\beta,\gamma)$ (Phase \textsf{1}), i.e. as links between individuals become strong enough, the ``neutral'' solution $\underline m_{*}^{0}$ becomes unstable (although it possesses a stable manifold), and one of the two opinions eventually prevails ({\em polarization of opinions}). This means that the two solutions in part 3. of Theorem \ref{solutions}, different from $\underline m_{*}^{0}$, are stable; the state space, besides the stable manifold for $\underline m_{*}^{0}$, gets partitioned into two parts, each attracted by one of the two stable solutions. This picture has been proved rigorously for the homogeneous case in \cite{DPRST09}, but it is  well supported by numerical evidence also in the inhomogeneous case. \\
Phase \textsf{2} is absent for $h=0$, thus it is a genuine effect of the inhomogeneity. If, from Phase \textsf{1}, we increase the link between individuals and their reference groups (specifically, as $h$ crosses the value $h(\beta,\gamma)$), then stability of the neutral solution is recovered. However, at least for moderate $h$, stability of one solution with $m^{\sigma}_* >0$ is maintained: in other words, polarization of opinions may occur or not, depending on the initial condition.

\subsection{Analysis of fluctuations: Central Limit
Theorem}

Thanks to Theorem \ref{lln} we established a law of large numbers for the empirical measure $\rho_N$, that is $\rho_N \longrightarrow Q_{*}$. We are going to analyze the Normal fluctuations around the limit $Q_{*}$. We are also interested in the $N$-asymptotic distribution of $\rho_N - Q_{*}$.\\
Using a weak convergence-type approach based on uniform convergence of the infinitesimal generators, deeply explained in \cite{EtKu86}, it is possible to provide a dynamical interpretation of the recalled law of large numbers.
 
Let $f:\mathscr{S}^2\longrightarrow\mathbb{R}$ be a function and define $\rho_N(t)$, the marginal distribution of $\rho_N$ at time $t$, by
\[\int f(\sigma, \omega) \, d\rho_N(t) = \frac{1}{N} \sum_{j=1}^{N} f (\sigma_j(t), \omega_j(t))\,.\]
We have $m^{\underline{\sigma}}_{N}(t) = m^{\underline{\sigma}}_{\rho_N(t)}$. For each fixed $t$, $\rho_N(t)$ is a probability on $\mathscr{S}^2$ and so, by the considerations which led as to introduce the expectations \eqref{mag2d'}, \eqref{mag2a'}, \eqref{mag2b'} and \eqref{mag2c'}, we can proceed similarly saying that $\rho_N(t)$ is completely determined by the vector $(m^{\underline{\eta}}_{\rho_N(t)}$, $m^{\underline{\sigma}}_{\rho_N(t)}$, $m^{\underline{\omega} }_{\rho_N(t)}$, $m^{\underline{\sigma} \, \underline{\omega}}_{\rho_N(t)}$,  $m^{\underline{\sigma} \, \underline{\eta}}_{\rho_N(t)}$, $m^{\underline{\omega} \, \underline{\eta}}_{\rho_N(t)}$, $m^{\underline{\sigma} \, \underline{\omega} \, \underline{\eta}}_{\rho_N(t)})$ and seeing it as a seven-dimensional object.\\ Thus $(\rho_N(t))_{t\in[0,T]}$ is a seven-dimensional flow. 

A simple consequence of Theorem \ref{lln} is the following convergence of flows:
\begin{equation}\label{flx'}
(\rho_N(t))_{t\in[0,T]} \longrightarrow (q_t)_{t\in[0,T]} \,,
\end{equation}
where the convergence is meant in probability, with respect to the weak topology for measure-valued processes. Since the flow of marginals contains less information than the full measure of paths, the law of large numbers in \eqref{flx'} is weaker than the one in Theorem \ref{lln}.\\ However, the corresponding fluctuation flow
\[ (N^{1/2}(\rho_N(t) - q_t))_{t\in[0,T]}\]
is also a finite-dimensional flow, whose limiting distribution can be explicitly determined.

\begin{TEO}\label{TLC}
In the limit $N\longrightarrow +\infty$, the seven-dimensional fluctuation process\\ $(r_N(t), x_N(t), y_N(t), z_N(t), u_N(t), v_N(t), w_N(t))$, defined by
\[ 
r_N(t) := N^{1/2} m^{\underline{\eta}}_{\rho_N(t)} 
\]
\begin{equation*}
\begin{split}
%\[ 
x_N(t)  &:=  N^{1/2} \left( m^{\underline{\sigma}}_{\rho_N(t)} - m^{\sigma}_{t} \right) \qquad u_N(t)  :=  N^{1/2} \left( m^{\underline{\sigma} \, \underline{\eta}}_{\rho_N(t)}- m^{\sigma \eta}_{t} \right)\\
%\]
%\[ 
y_N(t)  &:=  N^{1/2} \left( m^{\underline{\omega}}_{\rho_N(t)}- m^{\omega}_{t} \right) \qquad v_N(t)  :=  N^{1/2} \left( m^{\underline{\omega} \, \underline{\eta}}_{\rho_N(t)}- m^{\omega \eta}_{t} \right)\\
%\]
%\[
z_N(t)  &:=  N^{1/2} \left( m^{\underline{\sigma} \, \underline{\omega}}_{\rho_N(t)}- m^{\sigma \omega}_{t} \right) \quad\,\, w_N(t)  :=  N^{1/2} \left( m^{\underline{\sigma} \, \underline{\omega} \, \underline{\eta}}_{\rho_N(t)}- m^{\sigma \omega \eta}_{t} \right)
%\]
\end{split}
\end{equation*}
converges (in the sense of weak convergence of stochastic processes) to a limiting seven-dimensional Gaussian process $(r(t), x(t), y(t), z(t), u(t), v(t), w(t))$, which is the unique solution of the linear stochastic differential equation
\[
\begin{array}{ccl}
    dr(t) & = & 0 \\
    \begin{bmatrix}
        dx(t)\\
        dy(t)\\
        dz(t)\\
        du(t)\\
        dv(t)\\
        dw(t)
    \end{bmatrix}
    & = & 2 \, \mathscr{H} \, A_1(t)dt+2A_2(t)
    \begin{bmatrix}
        x(t)\\
        y(t)\\
        z(t)\\
        u(t)\\
        v(t)\\
        w(t)
    \end{bmatrix}
    dt+D(t)
    \begin{bmatrix}
        dB_1(t)\\
        dB_2(t)\\
        dB_3(t)\\
        dB_4(t)\\
        dB_5(t)\\
        dB_6(t)
    \end{bmatrix} 
\end{array}\,,
\]
where $B_1, B_2, B_3, B_4, B_5, B_6$ are independent Standard Brownian motions, and $\mathscr{H}$ is a Standard Gaussian random variable, % non devo dimostrare che H e i moti browniani sono indipendenti perche' 1. lo ottengo a posteriori per la teoria delle diffusioni 2. H dipende solo dalle condizioni iniziali
\begin{displaymath}
    A_1(t)=
      \begin{bmatrix}
      0\\
      \sinh(\gamma h)\cosh(\gamma m^{\sigma}_t) \\
      0\\
      0\\
      \cosh(\gamma h)\sinh(\gamma m^{\sigma}_t) \\
      \sinh(\beta)
      \end{bmatrix}\,,
\end{displaymath}
$A_2(t)=$
{\tiny
\[ %{\tiny
\begin{bmatrix}
-\mathrm{ch}(\beta) & \mathrm{sh}(\beta) & 0 & 0 & 0 & 0 \\ 
&&&&&\\
\begin{array}{c} -\gamma m^{\omega}_t \mathrm{ch}(\gamma h) \mathrm{sh}(\gamma m^{\sigma}_t) \\ -\gamma m^{\omega \eta}_t  \mathrm{sh}(\gamma h) \mathrm{ch}(\gamma m^{\sigma}_t) \\ + \gamma \mathrm{ch}(\gamma h) \mathrm{ch}(\gamma m^{\sigma}_t) \end{array} & -\mathrm{ch}(\gamma h) \mathrm{ch}(\gamma m^{\sigma}_t) & 0 & 0 & -\mathrm{sh}(\gamma h) \mathrm{sh}(\gamma m^{\sigma}_t) & 0 \\ 
&&&&&\\
\begin{array}{c} \mathrm{ch}(\gamma h) \mathrm{sh}(\gamma m^{\sigma}_t) \\ + \gamma m^{\sigma}_t \mathrm{ch}(\gamma h)   \mathrm{ch}(\gamma m^{\sigma}_t) \\ - \gamma m^{\sigma \omega}_t \mathrm{ch}(\gamma h)  \mathrm{sh}(\gamma m^{\sigma}_t) \\ + \gamma  m^{\sigma \eta}_t \mathrm{sh}(\gamma h) \mathrm{sh}(\gamma m^{\sigma}_t)\\ - \gamma m^{\sigma \omega \eta}_t\mathrm{sh}(\gamma h)  \mathrm{sh}(\gamma m^{\sigma}_t) \end{array} & 0 & \begin{array}{c} -\mathrm{ch}(\beta)\\ -\mathrm{ch}(\gamma h) \mathrm{ch}(\gamma m^{\sigma}_t) \end{array}  & \mathrm{sh}(\gamma h) \mathrm{ch}(\gamma m^{\sigma}_t) & 0 & -\mathrm{sh}(\gamma h) \mathrm{sh}(\gamma m^{\sigma}_t) \\ 
&&&&&\\
0 & 0 & 0 & -\mathrm{ch}(\beta) & \mathrm{sh}(\beta) & 0 \\ 
&&&&&\\
\begin{array}{c} \gamma \mathrm{sh}(\gamma h) \mathrm{sh}(\gamma m^{\sigma}_t) \\ - \gamma m^{\omega}_t \mathrm{sh}(\gamma h) \mathrm{ch}(\gamma m^{\sigma}_t) \\ -\gamma m^{\omega \eta}_t  \mathrm{ch}(\gamma h) \mathrm{sh}(\gamma m^{\sigma}_t) \end{array} & -\mathrm{sh}(\gamma h) \mathrm{sh}(\gamma m^{\sigma}_t) & 0 & 0 & -\mathrm{ch}(\gamma h) \mathrm{ch}(\gamma m^{\sigma}_t) & 0 \\ 
&&&&&\\
\begin{array}{c} \mathrm{sh}(\gamma h) \mathrm{ch}(\gamma m^{\sigma}_t) \\ + \gamma m^{\sigma}_t \mathrm{sh}(\gamma h)   \mathrm{sh}(\gamma m^{\sigma}_t) \\ - \gamma m^{\sigma \omega}_t \mathrm{sh}(\gamma h)  \mathrm{sh}(\gamma m^{\sigma}_t) \\ + \gamma  m^{\sigma \eta}_t \mathrm{ch}(\gamma h) \mathrm{ch}(\gamma m^{\sigma}_t)\\ - \gamma m^{\sigma \omega \eta}_t\mathrm{ch}(\gamma h)  \mathrm{sh}(\gamma m^{\sigma}_t) \end{array} & 0 & -\mathrm{sh}(\gamma h) \mathrm{sh}(\gamma m^{\sigma}_t) & \mathrm{ch}(\gamma h) \mathrm{sh}(\gamma m^{\sigma}_t) & 0 & \begin{array}{c} -\mathrm{ch}(\beta)\\ -\mathrm{ch}(\gamma h) \mathrm{ch}(\gamma m^{\sigma}_t) \end{array}
\end{bmatrix}%} 
\,,\]
}
$D(t)$ is a suitable $6 \times 6$ matrix and $(r(0), x(0), y(0), z(0), u(0), v(0), w(0))$ has a centered Gaussian distribution with covariance matrix
\[\!\!\!\! {\tiny
\begin{bmatrix}
1 & 0 & 0 & 0 & 0 & 0 & 0 \\
&&&&&&\\
0 & 1-(m_{\lambda}^{\sigma})^2 & m_{\lambda}^{\sigma \omega} -m_{\lambda}^{\sigma}m_{\lambda}^{\omega} & m_{\lambda}^{\omega}-m_{\lambda}^{\sigma}m_{\lambda}^{\sigma \omega} & -m_{\lambda}^{\sigma}m_{\lambda}^{\sigma \eta} & -m_{\lambda}^{\sigma}m_{\lambda}^{\omega \eta} & -m_{\lambda}^{\sigma}m_{\lambda}^{\sigma \omega \eta}  \\
&&&&&&\\
0 & m_{\lambda}^{\sigma \omega} -m_{\lambda}^{\sigma}m_{\lambda}^{\omega} & 1- (m_{\lambda}^{\omega})^2 & m_{\lambda}^{\sigma} -m_{\lambda}^{\omega}m_{\lambda}^{\sigma\omega} & -m_{\lambda}^{\omega}m_{\lambda}^{\sigma\eta} & -m_{\lambda}^{\omega}m_{\lambda}^{\omega \eta} & -m_{\lambda}^{\omega}m_{\lambda}^{\sigma\omega\eta}  \\
&&&&&&\\
0 & m_{\lambda}^{\omega} -m_{\lambda}^{\sigma}m_{\lambda}^{\sigma\omega} & m_{\lambda}^{\sigma} -m_{\lambda}^{\omega}m_{\lambda}^{\sigma\omega} & 1-(m_{\lambda}^{\sigma\omega})^2 & -m_{\lambda}^{\sigma \omega}m_{\lambda}^{\sigma\eta} & -m_{\lambda}^{\sigma\omega}m_{\lambda}^{\omega\eta} & -m_{\lambda}^{\sigma\omega}m_{\lambda}^{\sigma\omega\eta} \\
&&&&&&\\
0 & -m_{\lambda}^{\sigma}m_{\lambda}^{\sigma\eta} & -m_{\lambda}^{\omega}m_{\lambda}^{\sigma\eta} & -m_{\lambda}^{\sigma \omega}m_{\lambda}^{\sigma\eta} & 1-(m_{\lambda}^{\sigma\eta})^2 & -m_{\lambda}^{\sigma\eta}m_{\lambda}^{\omega \eta} & -m_{\lambda}^{\sigma\eta}m_{\lambda}^{\sigma\omega\eta}\\
&&&&&&\\
0 & -m_{\lambda}^{\sigma}m_{\lambda}^{\omega\eta} & -m_{\lambda}^{\omega}m_{\lambda}^{\omega\eta} & -m_{\lambda}^{\sigma\omega}m_{\lambda}^{\omega\eta} & -m_{\lambda}^{\sigma\eta}m_{\lambda}^{\omega\eta} & 1-(m_{\lambda}^{\omega\eta})^2 & -m_{\lambda}^{\omega\eta}m_{\lambda}^{\sigma\omega\eta}\\
&&&&&&\\
0 & -m_{\lambda}^{\sigma}m_{\lambda}^{\sigma\omega\eta} & -m_{\lambda}^{\omega}m_{\lambda}^{\sigma\omega\eta} & -m_{\lambda}^{\sigma\omega}m_{\lambda}^{\sigma\omega\eta} & -m_{\lambda}^{\sigma\eta}m_{\lambda}^{\sigma\omega\eta} & -m_{\lambda}^{\omega\eta}m_{\lambda}^{\sigma\omega\eta} & 1-(m_{\lambda}^{\sigma\omega\eta})^2
\end{bmatrix}} \,.
\]
\end{TEO}

\Proof We only sketch the main tools useful to prove the Theorem; we avoid the complete computations, since they are straightforward and very similar to those of the homogeneous case, which can be found in ~\cite{DPRST09}.

First, we have to prove that $(m^{\underline{\eta}}_{\rho_N(t)}$, $m^{\underline{\sigma}}_{\rho_N(t)}$, $m^{\underline{\omega}}_{\rho_N(t)}$, $m^{\underline{\sigma} \, \underline{\omega}}_{\rho_N(t)}$, $m^{\underline{\sigma} \, \underline{\eta}}_{\rho_N(t)}$, $m^{\underline{\omega} \, \underline{\eta}}_{\rho_N(t)}$, $m^{\underline{\sigma} \, \underline{\omega} \, \underline{\eta}}_{\rho_N(t)})$ is an order parameter for the model, i.e. its evolution is Markovian. So, let denote by $\mathcal{K}_N$ the infinitesimal generator of this process. We apply the operator \eqref{IG1'} to a function $\phi$, which is a composition of functions of the type: 
\[
\phi(m^{\underline{\eta}}_{\rho_N}, m^{\underline{\sigma}}_{\rho_N}, m^{\underline{\omega}}_{\rho_N}, m^{\underline{\sigma} \, \underline{\omega}}_{\rho_N}, m^{\underline{\sigma} \, \underline{\eta}}_{\rho_N}, m^{\underline{\omega} \, \underline{\eta}}_{\rho_N}, m^{\underline{\sigma} \, \underline{\omega} \, \underline{\eta}}_{\rho_N}):\mathscr{S}^{2N} \longrightarrow \mathbb{R},
\] 
where we are considering $(\underline{\sigma}, \underline{\omega})$ as a variable. The goal is we obtain the same applying another infinitesimal generator  (which is exactly the operator $\mathcal{K}_N$) to a function 
\[
\phi(m^{\underline{\eta}}_{\rho_N}, m^{\underline{\sigma}}_{\rho_N}, m^{\underline{\omega}}_{\rho_N}, m^{\underline{\sigma} \, \underline{\omega}}_{\rho_N}, m^{\underline{\sigma} \, \underline{\eta}}_{\rho_N}, m^{\underline{\omega} \, \underline{\eta}}_{\rho_N}, m^{\underline{\sigma} \, \underline{\omega} \, \underline{\eta}}_{\rho_N}): [-1,+1]^7 \longrightarrow \mathbb{R},
\] 
where $m^{\underline{\eta}}_{\rho_N}, m^{\underline{\sigma}}_{\rho_N}, m^{\underline{\omega}}_{\rho_N}, m^{\underline{\sigma} \, \underline{\omega}}_{\rho_N}, m^{\underline{\sigma} \, \underline{\eta}}_{\rho_N}, m^{\underline{\omega} \, \underline{\eta}}_{\rho_N}, m^{\underline{\sigma} \, \underline{\omega} \, \underline{\eta}}_{\rho_N}$ are seen as variables now.\\
In other words, it can be shown that
\begin{align*}
L_N(\phi(m^{\underline{\eta}}_{\rho_N}, m^{\underline{\sigma}}_{\rho_N}, m^{\underline{\omega}}_{\rho_N}, & m^{\underline{\sigma} \, \underline{\omega}}_{\rho_N}, m^{\underline{\sigma} \, \underline{\eta}}_{\rho_N}, m^{\underline{\omega} \, \underline{\eta}}_{\rho_N}, m^{\underline{\sigma} \, \underline{\omega} \, \underline{\eta}}_{\rho_N})) = \nonumber \\
& = (\mathcal{K}_N \phi)(m^{\underline{\eta}}_{\rho_N}, m^{\underline{\sigma}}_{\rho_N}, m^{\underline{\omega}}_{\rho_N}, m^{\underline{\sigma} \, \underline{\omega}}_{\rho_N}, m^{\underline{\sigma} \, \underline{\eta}}_{\rho_N}, m^{\underline{\omega} \, \underline{\eta}}_{\rho_N}, m^{\underline{\sigma} \, \underline{\omega} \, \underline{\eta}}_{\rho_N}) \nonumber \,.
\end{align*}
The seven-dimensional fluctuation process $(r_N(t), x_N(t), y_N(t), z_N(t), u_N(t), v_N(t), w_N(t))$ is a Markov process too, since it is a deterministic and invertible function of $(m^{\underline{\eta}}_{\rho_N}, m^{\underline{\sigma}}_{\rho_N}, m^{\underline{\omega}}_{\rho_N},$ $m^{\underline{\sigma} \, \underline{\omega}}_{\rho_N}, m^{\underline{\sigma} \, \underline{\eta}}_{\rho_N}, m^{\underline{\omega} \, \underline{\eta}}_{\rho_N}, m^{\underline{\sigma} \, \underline{\omega} \, \underline{\eta}}_{\rho_N})$. With the same reasoning as before, we can find the explicit expression of the infinitesimal generator $\mathcal{H}_N$, driving the dynamics of this fluctuation process.\\
Now, it can be found the limiting generator $\mathcal{H}$ of $\mathcal{H}_N$ and, using the results about the convergence of stochastic processes developed in ~\cite{EtKu86} (Chapter 4, Corollary 8.7), we can conclude our fluctuation process converges weakly to a Gaussian one, whose dynamics is driven by $\mathcal{H}$ and which solves the diffusion equation given in the statement of the Theorem.
\fine

\section{Critical dynamics ($\gamma = \cosh^2(\gamma h)/\tanh(\beta)$)}
We are going to consider the ``critical dynamics'' of the system, in other words the long-time behavior of the fluctuations in the threshold case, when $\gamma = \frac{\cosh^2(\gamma h)}{\tanh(\beta)}$. In the previous section we have seen that in a time interval $[0,T]$, where $T$ is fixed, and in the infinite volume limit, we have Normal fluctuations for the system. Indeed, the infinitesimal generator of the rescaled process converges to the infinitesimal generator of a diffusion and the rescaled process itself converges weakly to that diffusion. It means we can provide a Central Limit Theorem for all the values of $\beta$ and $\gamma$. This Central Limit Theorem continues to be valid in the critical case, but there is an eigenvalue of the covariance matrix $\Sigma_t$ which grows polynomially in $t$ and identifies the critical direction. This fact implies that the size of the Normal fluctuations must be further rescaled (in space and in time), because their size around the deterministic limit increases in time. The presence of the constant drift in the dynamics of the Normal fluctuations influences the construction and the behavior of the critical fluctuation process; in fact, it forces us to rescale the time by a smaller power of $N$ than it would be done in the homogeneous model.  The limiting process of these fluctuations is still Gaussian, since solution of a deterministic equation with constant (but random) drift given by a Gaussian random variable.  

First of all, we need to locate the critical direction in the seven-dimensional space of the order parameters. In the rest of the section, we will consider $\gamma = \frac{\cosh^2(\gamma h)}{\tanh(\beta)}$ and let us assume that the initial condition $\lambda$ is a product measure such
that
\[
m_t^{\sigma} = 0, \quad m_t^{\omega} = 0, \quad m_t^{\sigma \omega} =
\frac{\tanh(\beta) \tanh(\gamma h) \sinh(\gamma h) + \sinh(\beta)}{\cosh(\beta) +
\cosh(\gamma h)}
\]
\[
m_t^{\sigma \eta} = \tanh(\beta) \tanh(\gamma h), \quad m_t^{\omega \eta} =
\tanh(\gamma h), \quad m_t^{\sigma \omega \eta} = 0\,,
\]
for every value of $t\geq 0$. Observe that if it holds at initial time, it will be
true for all subsequent times by stationarity. \\

In the critical case, when $\gamma = \frac{\cosh^2(\gamma h)}{\tanh(\beta)}$, the matrix $\Sigma_{t}$ has an eigenvalue growing polynomially in $t$ . The critical direction is determined by the right eigenvector corresponding to the eigenvalue increasing to infinity of $\Sigma_{t}$, which is also the right eigenvector corresponding to the null eigenvalue of $A_2$, the drift matrix of Theorem \ref{TLC} (independent of $t$ under our assumptions).
This matrix can not be completely diagonalized, but it is possible to find a basis of generalized right eigenvectors $\{\underline{a}_1$,\dots, $\underline{a}_6\}$, allowing us to reduce it in Jordan canonical form. Let $A$ be the matrix whose rows are $\underline{a}_1$,\dots, $\underline{a}_6$. It is convenient to consider the following \emph{change of variables} (it will be used to construct the critical fluctuation process in Theorem \ref{thmCRTDYN'}):
\[
\begin{array}{ccl}
    r_N(t) & = & r_N(t) \\
    \begin{bmatrix}
        \bar{x}_N(t)\\
        \bar{y}_N(t)\\
        \bar{z}_N(t)\\
        \bar{u}_N(t)\\
        \bar{v}_N(t)\\
        \bar{w}_N(t)
    \end{bmatrix}
    & = & 
    \begin{bmatrix}
    N^{-1/4} & 0 & 0 & 0 & 0 & 0 \\
    0 & N^{-1/4} & 0 & 0 & 0 & 0 \\
    0 & 0 & N^{-1/4} & 0 & 0 & 0 \\
    0 & 0 & 0 & N^{-1/4} & 0 & 0 \\
    0 & 0 & 0 & 0 & N^{-1/4} & 0 \\
    0 & 0 & 0 & 0 & 0 & N^{-1/4} \\
    \end{bmatrix}
    A
    \begin{bmatrix}
        x_N(t)\\
        y_N(t)\\
        z_N(t)\\
        u_N(t)\\
        v_N(t)\\
        w_N(t)
    \end{bmatrix} 
\end{array}\,.
\]
If we set $\underline{a}_1$ to be the eigenvector corresponding to the null eigenvalue, by this analysis we obtain the following critical direction:
\begin{align*}
\bar{x}_N(t) &=  N^{-1/4} \; \underline{a}_1 \cdot (x_N(t),y_N(t),z_N(t),u_N(t),v_N(t),w_N(t))' \\
&= N^{1/4} \Big[\cosh(\gamma h) m^{\underline{\sigma}}_{\rho_N(t)} + \sinh(\beta) m^{\underline{\omega}}_{\rho_N(t)} \Big]\,.
\end{align*}

\begin{REM}
Notice that the critical direction $\bar{x}$ does not depend on the random environment and it is one-dimensional.
\end{REM}  

\begin{DEF}
We say that the sequence of stochastic processes $\{X_n(t)\}_{n}$, for $t \in [0,T]$, \emph{collapses to zero} if for every $\varepsilon >0$, 
\[
\lim_{n \to +\infty} P \left( \, \sup_{t \in [0,T]} \vert X_n(t) \vert > \varepsilon \right)=0 \,.
\]
\end{DEF}

\begin{TEO}\label{thmCRTDYN'}
For $t \in [0,T]$, let consider the seven-dimensional critical fluctuation process
\begin{allowdisplaybreaks}
\begin{align}\label{fluctcrt}
%\begin{split}
r_{N}(t) &= N^{1/2} m^{\underline{\eta}}_{\rho_N(t)} \nonumber\\
\bar{x}_{N}(t) &= N^{1/4} \Big[\cosh(\gamma h) m^{\underline{\sigma}}_{\rho_N(t)} + \sinh(\beta) m^{\underline{\omega}}_{\rho_N(t)} \Big] \nonumber\\
\bar{y}_{N}(t) &= N^{1/4} \Big[(\cosh(\gamma h) - \cosh(\beta)) m^{\underline{\sigma} \, \underline{\eta}}_{\rho_N(t)} + \sinh(\beta) m^{\underline{\omega} \, \underline{\eta}}_{\rho_N(t)} \nonumber \\
      & \qquad\qquad - (\cosh(\gamma h) - \cosh(\beta)) \tanh(\beta) \tanh(\gamma h) - \sinh(\beta)\tanh(\gamma h) \Big] \nonumber \\
\bar{z}_{N}(t) &= N^{1/4} \Big[ m^{\underline{\omega} \, \underline{\eta}}_{\rho_N(t)} - \tanh(\gamma h)\Big] \nonumber\\
\bar{u}_{N}(t) &= N^{1/4} \Bigg[ m^{\underline{\sigma} \, \underline{\omega}}_{\rho_N(t)} - \tanh(\gamma h) m^{\underline{\sigma} \, \underline{\eta}}_{\rho_N(t)} + \tanh(\beta) \tanh(\gamma h) m^{\underline{\omega} \, \underline{\eta}}_{\rho_N(t)} \nonumber \\
      & \hspace{5cm} - \frac{\tanh(\beta) \tanh(\gamma h) \sinh(\gamma h) + \sinh(\beta)}{\cosh(\beta) + \cosh(\gamma h)}\Bigg] \\
%\end{split}
%\end{equation}
%\begin{equation*}
%\begin{split}              
\bar{v}_{N}(t) &= N^{1/4} \Big[\, 2 \cosh(\beta) \sinh(\gamma h) (\cosh(\beta) + 2 \cosh(\gamma h)) m^{\underline{\sigma}}_{\rho_N(t)} \nonumber \\
       & \hspace{3.3cm} - 2 \sinh(\beta) \sinh(\gamma h) (\cosh(\beta) + 2 \cosh(\gamma h)) m^{\underline{\omega}}_{\rho_N(t)} \Big] \nonumber \\
       %\end{split}
%\end{equation}
%\begin{equation*}
%\begin{split} 
\bar{w}_{N}(t) &= N^{1/4} \Big[ - \tanh(\beta) \sinh(\gamma h) (\cosh(\beta) + 2 \cosh(\gamma h)) m^{\underline{\omega}}_{\rho_N(t)} \nonumber\\
       & \hspace{6.3cm} + (\cosh(\beta) + \cosh(\gamma h))^2 \, m^{\underline{\sigma} \, \underline{\omega} \, \underline{\eta}}_{\rho_N(t)}\Big]\,. \nonumber 
%\end{split}
\end{align}
\end{allowdisplaybreaks}\\
Then, as $N \longrightarrow +\infty$, $r_{N}(t)$ converges to $\mathscr{H}$, a Standard Gaussian random variable, the processes $\bar{y}_{N}(N^{1/4}t)$, $\bar{z}_{N}(N^{1/4}t)$, $\bar{u}_{N}(N^{1/4}t)$, $\bar{v}_{N}(N^{1/4}t)$, $\bar{w}_{N}(N^{1/4}t)$ collapse to zero and $\bar{x}_{N}(N^{1/4}t)$ converges, in the sense of weak convergence of stochastic processes, to a limiting Gaussian process 
\[
\bar{x}(N^{1/4}t) = 2 \, \mathscr{H} \sinh(\beta) \sinh(\gamma h)  \,  t \,.
\]
\end{TEO}

\Proof See Section \ref{proofs}. \fine

\begin{REM}
In view of Theorem \ref{thmCRTDYN'} is now clear why it is convenient the change of variables we introduced. In the new system of coordinates only one process survives the critical space-time scaling.
\end{REM}

The previous Theorem is valid even in the homogeneous model, when $h=0$, but the result is trivial: the limiting process $\bar{x}(N^{1/4}t) \equiv 0$. In that case, the absence of the constant drift, due to the random environment, allows us to amplify the time by a factor $N^{1/2}$ in order to have an appreciable and different result: the limiting critical fluctuation process is no more trivial and, besides, it is non-Gaussian. We state here this result, for the complete proof we refer to ~\cite{Sar07}.

\begin{TEO}\label{critico om}
Assume $h=0$. For $t \in [0,T]$, if we consider the critical fluctuation process
\begin{align}\label{fluctcrth'}
        \xi_N(t)  &:=  N^{1/4} \left[m^{\underline{\sigma}}_{\rho_N(N^{1/2}t)} - \tanh(\beta)m^{\underline{\omega}}_{\rho_N(N^{1/2}t)} \right] \nonumber \\
        \vartheta_N(t)  &:=  N^{1/4} \left[m^{\underline{\sigma}}_{\rho_N(N^{1/2}t)} + \sinh(\beta)m^{\underline{\omega}}_{\rho_N(N^{1/2}t)} \right]\\
        \zeta_N(t)  &:=  N^{1/4} \left[m^{\underline{\sigma\omega}}_{\rho_N(N^{1/2}t)} - \frac{\sinh(\beta)}{\cosh(\beta)+1} \right] \nonumber
\end{align}
then, as $N \longrightarrow +\infty$, $\xi_N(t),\zeta_N(t)$ collapse to zero and $\vartheta_N(t)$ converges, in the sense of weak convergence of stochastic processes, to a limiting non-Gaussian process $\vartheta(t)$, which is the unique solution of the following stochastic differential equation:
\[ \left\{
\begin{array}{ccl}
 d\vartheta(t)  &=&   -\displaystyle{\frac{2\cosh^3(\beta)}{3 \sinh^2(\beta)[\cosh(\beta)+1]^3}} \, \vartheta^{3}(t) \, dt + 2 \cosh(\beta) \, dB(t) \\
 &&\\
 \vartheta(0) &=&0
\end{array} \right.
\]
where $B$ is a standard Brownian motion.
\end{TEO}

\section{Conclusions}

We have proposed a model for social interactions having the following main features:
\begin{itemize}
\item[$\RHD$]
Interaction is of mean field type: the same information is available to all individuals.
\item[$\RHD$]
Individuals are not identical (inhomogeneity): they are divided into reference groups. 
\item[$\RHD$]
Equilibrium dynamics are not time reversible.
\end{itemize}
We have shown various asymptotic results in the limit as the number $N$ of individuals goes to infinity. First we have shown a law of large numbers, that describes the dynamics of a specific individual in the limit of an infinite community ($N \ra +\infty$). The large time behavior of this dynamics exhibit phase transitions: depending on the parameters of the model, and possibly on the initial condition, states of individuals may either polarize - due to strong tendency to conformism - or tend to a ``neutral'' configuration, where states are mainly dictated by the reference group individuals belong to, and only weakly influenced by interactions within the community. We have then proved a Central Limit Theorem, which provides explicit normal corrections to the $N \ra +\infty$ limiting dynamics. Finally, fluctuations around the limiting dynamics have been studied in more details in the critical region of parameters, separating the phases of polarization and non-polarization. In this region, large-time fluctuations exhibit peculiar scaling properties, which are strongly influenced by  the inhomogeneity of the community.

\section{Proofs}\label{proofs}

\subsection{Proof of Lemma \ref{MKV2'}}

By definition \eqref{mag2a'}, Proposition \ref{thmMKV'} and Theorem \ref{lln} we deduce that
\begin{allowdisplaybreaks}
\begin{align*}
\dot{m}^{\sigma}_{t} \! = & \! \sum_{\sigma, \omega \in \mathscr{S}} \sigma\, \dot{q}_t(\sigma, \omega) = \frac{1}{2} \sum_{\sigma, \omega \in \mathscr{S}} \sum_{\eta \in \mathscr{S}} \sigma \, \dot{q}_t^{\eta} (\sigma, \omega) = \frac{1}{2} \sum_{\sigma, \omega \in \mathscr{S}} \sum_{\eta \in \mathscr{S}} \sigma \,\mathcal{L}^{\eta} q_t^{\eta} (\sigma, \omega) \\
                     = & \frac{1}{2} \! \sum_{\sigma, \omega \in \mathscr{S}} \sum_{\eta \in \mathscr{S}} \sigma \left\{ \nabla^\sigma \left[ e^{-\beta \sigma \omega} q_t^{\eta} (\sigma, \omega)\right] + \nabla^\omega \left[ e^{-\gamma \omega ( m^{\sigma}_{t}+ h \eta)} q_t^{\eta} (\sigma, \omega) \right] \right\}\\
                     = & \frac{1}{2} \! \sum_{\sigma, \omega \in \mathscr{S}} \sum_{\eta \in \mathscr{S}} \sigma \left\{ e^{\beta \sigma \omega} q_t^{\eta}(-\sigma, \omega) -  e^{-\beta \sigma\omega} q_t^{\eta}(\sigma, \omega) + \nabla^\omega \left[ e^{-\gamma \omega ( m^{\sigma}_{t}+ h \eta)} q_t^{\eta} (\sigma, \omega) \right] \right\}\\
                     = & - \! \sum_{\sigma, \omega \in \mathscr{S}} \sum_{\eta \in \mathscr{S}} \sigma e^{- \beta \sigma \omega} q_t^{\eta} (\sigma, \omega) + \frac{1}{2}  \underbrace{\sum_{\sigma, \omega \in \mathscr{S}} \sum_{\eta \in \mathscr{S}} \sigma  \nabla^\omega \left[ e^{-\gamma \omega ( m^{\sigma}_{t}+ h \eta)} q_t^{\eta} (\sigma, \omega) \right]}_{= 0} \\
                     = & - \! \sum_{\sigma, \omega \in \mathscr{S}} \sum_{\eta \in \mathscr{S}}  \sigma \left[ \cosh(\beta) -\sigma \omega \sinh(\beta) \right]  q_t^{\eta}(\sigma, \omega) \\
                     = & - \cosh(\beta) \sum_{\sigma, \omega \in \mathscr{S}} \sum_{\eta\in\mathscr{S}} \sigma \, q_t^{\eta}(\sigma, \omega) + \sinh(\beta) \sum_{\sigma, \omega \in \mathscr{S}} \sum_{\eta\in\mathscr{S}} \omega \, q_t^{\eta}(\sigma, \omega)  \\
                     = &  -2 \, m^{\sigma}_{t} \cosh(\beta)  + 2\, m^{\omega}_{t} \sinh(\beta)\,,
\end{align*}
\end{allowdisplaybreaks}\\
where the last equality holds thanks to \eqref{mag2a'} and \eqref{mag2b'}. So the first equation of  Lemma \ref{MKV2'} is proved. Similarly, we can obtain all the others.

\subsection{Proof of Theorem \ref{solutions}}

The fact that $\underline{m}_*^0$ is an equilibrium point for all values of the parameters is easily shown by equations (\ref{soleq'}).
Now we are going to study the linear stability of this equilibrium. Denoting by
\[
\begin{array}{rccc} 
V: & [-1,1]^6 & \longrightarrow & \mathbb{R}^6 \\ 
    & \underline{x} := (x_1, x_2, x_3, x_4, x_5, x_6) & \longmapsto & (V_1(\underline{x}), V_2(\underline{x}), V_3(\underline{x}), V_4(\underline{x}), V_5(\underline{x}), V_6(\underline{x}))\,,
\end{array}
\] 
with
\begin{allowdisplaybreaks}
\begin{align*}
V_1(\underline{x}) & := -2\, x_1 \cosh(\beta ) + 2 \, x_2 \sinh(\beta )\\
V_2(\underline{x}) & := -2 \, x_2 \cosh(\gamma h) \cosh(\gamma x_1) - 2 \, x_5 \sinh(\gamma h) \sinh(\gamma x_1)\\
                           & \qquad + 2 \, \cosh(\gamma h) \sinh(\gamma x_1)\\
V_3(\underline{x}) &:= 2 \, x_1 \cosh(\gamma h) \sinh(\gamma x_1) - 2 \, x_3 \left[ \cosh(\beta) + \cosh(\gamma h) \cosh(\gamma x_1) \right]\\
                           & \qquad + 2 \, x_4 \sinh(\gamma h) \cosh(\gamma x_1) - 2 \, x_6 \sinh(\gamma h) \sinh(\gamma x_1)+2 \sinh(\beta)\\
V_4(\underline{x}) &:=  -2\, x_4 \cosh(\beta) + 2\, x_5 \sinh(\beta) \\
V_5(\underline{x}) &:= -2 \, x_2 \sinh(\gamma h) \sinh(\gamma x_1) - 2 \, x_5 \cosh(\gamma h) \cosh(\gamma x_1)\\
                           & \qquad  + 2 \, \sinh(\gamma h) \cosh(\gamma x_1) \\
V_6(\underline{x}) &:= 2 \, x_1 \sinh(\gamma h) \cosh(\gamma x_1) - 2 \, x_3 \sinh(\gamma h) \sinh(\gamma x_1)\\
                           & \qquad + 2 \, x_4 \cosh(\gamma h) \sinh(\gamma x_1) - 2 \, x_6 \left[ \cosh(\beta) + \cosh(\gamma h) \cosh(\gamma x_1) \right]\,,
\end{align*}
\end{allowdisplaybreaks}\\
the vector field of the system in Lemma \ref{MKV2'}, we obtain the linearized matrix evaluated in the stationary solution is $DV(\underline m_{*}^{0})$:
{\scriptsize
 \[
DV(\underline m_{*}^{0}) = 2 %{\scriptsize 
\begin{bmatrix}
-\mathrm{ch}(\beta) & \mathrm{sh}(\beta) & 0 & 0 & 0 & 0 \\
\frac{ \gamma}{\mathrm{ch}(\gamma h)}  & - \mathrm{ch}(\gamma h) & 0 & 0 & 0 & 0 \\
0 & 0 & - [\mathrm{ch}(\beta) + \mathrm{ch}(\gamma h)] &  \mathrm{sh}(\gamma h) & 0 & 0 \\
0 & 0 & 0 & - \mathrm{ch}(\beta) & \mathrm{sh}(\beta) & 0 \\
0 & 0 & 0 & 0 & -\mathrm{ch}(\gamma h) & 0 \\
\mathrm{sh}(\gamma h) + \gamma \frac{\mathrm{th}(\beta) \mathrm{th}(\gamma h)}{\mathrm{ch}(\beta) + \mathrm{ch}(\gamma h)} & 0 & 0 & 0 & 0 & - [\mathrm{ch}(\beta) + \mathrm{ch}(\gamma h)]\,.
\end{bmatrix}.%}
\]
}

Its eigenvalues are given by
\begin{align*}
\lambda_1 =& - \cosh(\beta) - \cosh(\gamma h) + \sqrt{[\cosh(\beta) - \cosh(\gamma h)]^2 + 4 \gamma \frac{\sinh(\beta)}{\cosh(\gamma h)}} \\
\lambda_2 =& - \cosh(\beta) - \cosh(\gamma h) - \sqrt{[\cosh(\beta) - \cosh(\gamma h)]^2 + 4 \gamma \frac{\sinh(\beta)}{\cosh(\gamma h)}} \\
\lambda_3 =& \lambda_4 = -2 [\cosh(\beta) + \cosh(\gamma h)] \\
\lambda_5 =& -2 \cosh(\gamma h) \\ 
\lambda_6 =& -2 \cosh(\beta)\,.
\end{align*} 
They are all real and it is easy to see that $\lambda_2, \lambda_3, \lambda_4, \lambda_5, \lambda_6<0$, for every value of $\beta$, $\gamma$, $h$; instead, the value of $\lambda_1$ depends on the parameters:
\begin{itemize}
\item[$\RHD$] if $\gamma < \frac{\cosh^2(\gamma h)}{\tanh(\beta)}$, then $\lambda_1<0$ and thus $\underline m_{*}^{0}$ is linearly stable;
\item[$\RHD$] if $\gamma = \frac{\cosh^2(\gamma h)}{\tanh(\beta)}$, then $\lambda_1=0$ and thus $DV(\underline m_{*}^{0})$ has a neutral direction;
\item[$\RHD$] if $\gamma > \frac{\cosh^2(\gamma h)}{\tanh(\beta)}$, then $\lambda_1>0$ and thus the linearized system admits a direction which is unstable.
\end{itemize}

Having established linear stability of $\underline{m}_*^0$, we now look for further equilibria.
To this purpose, it is sufficient to study the behavior of the self-consistency relation satisfied by $m_{*}^{\sigma}$. Looking at the first expression in \eqref{soleq'}, we can write
\begin{equation}\label{selfcon'}
m_{*}^{\sigma} =  \Gamma_{\beta, \gamma, h} (m_{*}^{\sigma}) \,,
\end{equation}
where $\Gamma_{\beta,\gamma,h}(m_*^\sigma)$ is defined by \eqref{Gamma}. 
It follows from \eqref{selfcon'} that
\begin{itemize}
\item[$\RHD$] $m^{\sigma}_{*} \longmapsto \Gamma_{\beta, \gamma, h} (m_{*}^{\sigma})$ is a continuous function for all the values of $\beta$, $\gamma$ and $h$;
\item[$\RHD$] $\displaystyle{\lim_{m_{*}^{\sigma} \rightarrow \pm \infty} \Gamma_{\beta, \gamma, h} (m_{*}^{\sigma}) = \pm \tanh(\beta)};$
\item[$\RHD$] $\displaystyle{\Gamma^{\prime}_{\beta, \gamma, h} (m_{*}^{\sigma}) = \gamma \tanh(\beta) \frac{[1 + 2 \sinh^2(\gamma h)] \cosh^2(\gamma m_{*}^{\sigma}) - \sinh^2(\gamma h)}{[\cosh^2(\gamma m_{*}^{\sigma}) + \sinh^2(\gamma h)]^2} > 0}\,\,$, for every $\beta$, $\gamma$ and $h$.
\end{itemize}
Since $\Gamma_{\beta, \gamma, h} (m_{*}^{\sigma})$ is an odd function with respect to $m_{*}^{\sigma}$, we have $\Gamma_{\beta, \gamma, h} (0) = 0$ for all  $\beta$, $\gamma$ and $h$, so that \eqref{selfcon'} has the paramagnetic solution $m_{*}^{\sigma} = 0$ always. Now, we investigate under what conditions ferromagnetic solutions $m_{*}^{\sigma} > 0$ may occur. We restrict to work in the positive half-plane.\\
If
\begin{equation}\label{phtcond'}
\Gamma^{\prime}_{\beta, \gamma, h} (0) = \gamma  \frac{\tanh(\beta)}{\cosh^2(\gamma h)} > 1\,,
\end{equation}
then there is at least one ferromagnetic  solution. However, since $\Gamma_{\beta, \gamma, h} (m_{*}^{\sigma})$ is not always concave, there may be a ferromagnetic solution even when \eqref{phtcond'} fails. In this case, there must be at least two ferromagnetic solutions (corresponding to the curve $m^{\sigma}_{*} \longmapsto \Gamma_{\beta, \gamma, h} (m_{*}^{\sigma})$, crossing the diagonal first from below and then from above).\\
The regime defined by \eqref{phtcond'} lies under the curve \eqref{curve'}. An idea of when two ferromagnetic solutions arise may be obtained from the Taylor expansion of $\Gamma_{\beta, \gamma, h} (m_{*}^{\sigma})$ for small $m_{*}^{\sigma}$; in fact,
\[
\Gamma_{\beta, \gamma, h} (m_{*}^{\sigma}) = \gamma \frac{\tanh(\beta)}{\cosh^2(\gamma h)} m_{*}^{\sigma} + \gamma^3 \frac{\tanh(\beta)[2\cosh^2(\gamma h) - 3]}{3 \cosh^4(\gamma h)} (m_{*}^{\sigma})^3 + O \left( (m_{*}^{\sigma})^5 \right)
\]
and, on the curve defined by \eqref{curve'}, it reduces to
\[
\Gamma_{\beta, \gamma, h} (m_{*}^{\sigma}) =  m_{*}^{\sigma} + \gamma \left( \frac{2}{3} \gamma - \frac{1}{\tanh(\beta)} \right) (m_{*}^{\sigma})^3 + O \left( (m_{*}^{\sigma})^5 \right)\,,
\]
from which we can see that $\overline{\gamma} = \frac{3}{2 \tanh(\beta)}$ is a critical value. Indeed, if $\gamma > \overline{\gamma}$, then as $h$ increases through $h(\beta, \gamma)$ (i.e. $\Gamma^{\prime}_{\beta, \gamma, h} (0)$ decreases through 1), at least two positive ferromagnetic solutions occur, because $m^{\sigma}_{*} \longmapsto \Gamma_{\beta, \gamma, h} (m_{*}^{\sigma})$ is convex for small $m_{*}^{\sigma}$.

\subsection{Proof of Theorem \ref{thmCRTDYN'}}

%%%%%%%%%%%%%%%STEP 1%%%%%%%%%%%%%%%%%%%%%%%
\textbf{$\quad\;\,$STEP 1.} Let us denote by $\{\tau_{N}^M\}_{N \geq 1}$ a family of stopping times, defined as
\begin{align*}
\tau_{N}^M:=\inf_{t \geq 0} \, \{ &\, \vert \widetilde{x}_{N}(t) \vert \geq M \; \mathrm{or} \;\; \vert \widetilde{y}_{N}(t) \vert \geq M \; \mathrm{or} \;\; \vert \widetilde{z}_{N}(t) \vert \geq M \\
                                               &  \; \mathrm{or} \;\; \vert \widetilde{u}_{N}(t) \vert \geq M \; \mathrm{or} \;\; \vert \widetilde{v}_{N}(t) \vert \geq M \; \mathrm{or} \;\; \vert \widetilde{w}_{N}(t) \vert \geq M \}\,,
\end{align*}
where $M$ is a positive constant. We are interested in introducing such a sequence of stopping times, because in this way the processes $\widetilde{x}_N(t)$, $\widetilde{y}_N(t)$, $\widetilde{z}_N$, $\widetilde{u}_N(t)$, $\widetilde{v}_N(t)$, $\widetilde{w}_N(t)$  result to be bounded in the time interval $[0, T \wedge \tau_N^M]$; thanks to the Central Limit Theorem $r_N(t)$ is still bounded, in fact for every $\varepsilon>0$ there exists $M>0$ such that $ P \left\{ \vert r_N(t) \vert \geq M \right\} \leq \varepsilon $ for $t \in [0,T]$. \\

%%%%%%%%%%%%%%%%%%STEP 2%%%%%%%%%%%%%%%%%%%%%%%%%%
\textbf{STEP 2.} Now we prove that, for $t \in [0, T \wedge \tau_{N}^M]$, the non-critical directions converge to zero in probability and this implies that $\widetilde{y}_{N}(t), \widetilde{z}_{N}(t), \widetilde{u}_{N}(t), \widetilde{v}_{N}(t), \widetilde{w}_{N}(t) \longrightarrow 0$, as $N \longrightarrow \infty$. We show it only for the process $\widetilde{z}_{N}(t)$, because the calculations are analogous in the other cases. 
First we need the following technical results.

\begin{PROP}\label{lemmaMP}
For $t \in [0,T \wedge \tau_{N}^M]$, %if we consider only the space scaling
the process $(r_N(t)$, $\bar{x}_{N}(t)$, $\bar{y}_{N}(t)$, $\bar{z}_{N}(t)$, $\bar{u}_{N}(t)$, $\bar{v}_{N}(t)$, $\bar{w}_{N}(t))$, defined in \eqref{fluctcrt}, is a Markov process.
\end{PROP}

To prove that $(r_N(t), \bar{x}_{N}(t), \bar{y}_{N}(t), \bar{z}_{N}(t), \bar{u}_{N}(t), \bar{v}_{N}(t), \bar{w}_{N}(t) )$ is a Markov process, one must write down the expression of the infinitesimal generator $\mathcal{G}_N$, whose dynamics are driven by. To do it, we need the following technical Lemma.

\begin{LEMMA}\label{LmmGf}
Let $(X_t)_{t \geq 0}$ be a continuous-time Markov chain on a finite state space S, admitting an infinitesimal generator $L$. Let $g:S\longrightarrow S'$ be a given function, where $S'$ is a finite set. Assume that for every $f:S\longrightarrow\mathbb{R}$, $L(f \circ g)$ is a function of $g(x)$, i.e. $L(f \circ g)= (Kf) \circ g$. Then this last identity defines a linear operator $K$; moreover, $g(X_t)$ is a Markov process with infinitesimal generator $K$.
\end{LEMMA}

\Proof
Obviously $K$ is linear. Observing that
\begin{equation}\label{i}
e^{tL}(f \circ g) = (e^{tK}f) \circ g \,,
\end{equation}
we can conclude. In fact, $X_t$ is a Markov process with generator $L$, then we have
\begin{align*}
E[(f \circ g)(X_t) \vert X_0=x]       & \; = \; e^{tL} (f \circ g)(x) \\
                                                  &\stackrel{\mbox{{\tiny{\eqref{i}}}}}{=} e^{tK}f(g(x)) \\
                                                  & \;\, = \; E[f(g(X_t)) \vert g(X_0) = g(x)] 
\end{align*}
and the last inequality holds since $e^{tK}$ is a Markov semigroup and $L(f \circ g)= (Kf) \circ g$. Hence, $g(X_t)$ is a Markov process with infinitesimal generator $K$.
\fine

\emph{ \sl Proof of Proposition \ref{lemmaMP}. }
We apply Lemma \ref{LmmGf}. The process $\{(\underline{\sigma}(t), \underline{\omega}(t))\}_{t \geq 0}$ is a continuous-time Markov chain on the finite state space $\mathscr{S}^{2N}$, with infinitesimal generator $L_N$, defined by \eqref{IG1'}. 
Let consider the function
\[
\begin{array}{cccccc}
\zeta: & \mathscr{S}^{2N}       & \xrightarrow{\;\zeta_{1}\;} & [-1,+1]^7 & \xrightarrow{\;\zeta_{2}\;} & \mathbb{R}^7 \\
          & (\underline{\sigma},\underline{\omega}) & \longmapsto     & (m^{\underline{\eta}}_{\rho_N}, m^{\underline{\sigma}}_{\rho_N}, \dots, m^{\underline{\sigma} \, \underline{\omega} \, \underline{\eta}}_{\rho_N}) & \longmapsto & (r_N(t), \bar{x}_{N}(t), \dots, \bar{w}_{N}(t) )\,;
\end{array}
\]
it plays the role of $g$ in Lemma \ref{LmmGf}. Then, for every $\psi: \mathscr{S}^{2N} \longrightarrow \mathbb{R}$, 
we have
\[
L_N(\psi \circ \zeta) = (\mathcal{G}_N \psi) \circ \zeta
\]
and $\zeta(\underline{\sigma},\underline{\omega})$ is a Markov process with generator $\mathcal{G}_N$ given by
{\small{
\begin{allowdisplaybreaks}
\begin{align}\label{IG5'}
&\mathcal{G}_N \psi(r, \bar{x}, \bar{y}, \bar{z}, \bar{u}, \bar{v}, \bar{w}) = \!\!\!\!  \sum_{i, j, k\in\mathscr{S}} \vert A_{N}(i, j, k) \vert  e^{-\beta ij}\cdot\nonumber\\ 
&\,\cdot\bigg[\psi\bigg( r, \bar{x} - i\frac{2}{N^{3/4}} \mathrm{ch}(\gamma h),
\bar{y} - ik\frac{2}{N^{3/4}}(\mathrm{ch}(\gamma h) - \mathrm{ch}(\beta)),
\bar{z}, \bar{u} - ij \frac{2}{N^{3/4}} + ik\frac{2}{N^{3/4}}\mathrm{th}(\gamma h),  \nonumber\\
&%\hspace{75pt}
\quad\,\,\bar{v} - i\frac{4}{N^{3/4}}\mathrm{sh}(\gamma h) \mathrm{ch}(\beta) (\bb),\bar{w} - ijk \frac{2}{N^{3/4}} (\aa)^2 \bigg)  - \psi( r, \bar{x}, \bar{y},  \bar{z}, \bar{u}, \bar{v}, \bar{w})\bigg]\nonumber\\
&+ \!\!\!\! \sum_{i, j, k\in\mathscr{S}} \vert A_{N}(i, j, k) \vert e^{\gjm}\, \cdot  \nonumber\\ 
&\,\cdot \bigg[\psi\bigg( r,  \bar{x} - j\frac{2}{N^{3/4}}, \bar{y} - jk\frac{2}{N^{3/4}}\mathrm{sh}(\beta), \bar{z} - jk \frac{2}{N^{3/4}}, \bar{u} - ij\frac{2}{N^{3/4}}- jk\frac{2}{N^{3/4}}\mathrm{th}(\gamma h)\mathrm{th}(\beta),\nonumber\\ 
&\quad\,\,\bar{v} + j \frac{4}{N^{3/4}}\mathrm{sh}(\gamma h)\mathrm{sh}(\beta)(\bb),\bar{w} + j\frac{2}{N^{3/4}}\mathrm{sh}(\gamma h)\mathrm{th}(\beta)(\bb)\nonumber\\ 
&\hspace{200pt}
- ijk \frac{2}{N^{3/4}}(\aa)^2 \bigg) %- \nonumber\\
%&\hspace{240pt}
-\psi( r, \bar{x}, \bar{y},  \bar{z}, \bar{u}, \bar{v}, \bar{w}) \bigg] \,,
\end{align}
\end{allowdisplaybreaks}\\
}}
where $A_{N}(i,j,k)$ is the set of all triples $(\sigma_d, \omega_d, \eta_d)$, $d\in\{1,\dots,N\}$, such that $\sigma_d=i$, $\omega_d=j$ and $\eta_d=k$, with $i, j, k \in \mathscr{S}$; hence
\begin{allowdisplaybreaks}
{\small{
\begin{multline}\label{contatore}
\vert  A_N(i, j, k) \vert = \frac{N}{8}\bigg[ 
1 + \left(jk\mathrm{th}(\gamma h) + ik\mathrm{th}(\beta) \mathrm{th}(\gamma h) 
 + ij\frac{\mathrm{th}(\beta) \mathrm{th}(\gamma h) \mathrm{sh}(\gamma h) + \mathrm{sh}(\beta)}{\aa}\right)\\
 + k \frac{r}{N^{1/2}} 
   + \frac{\bar{x}}{N^{1/4} (\aa)}\left(i + j  \frac{\mathrm{ch}(\beta)}{\mathrm{sh}(\beta)} 
+ ijk \frac{\mathrm{sh}(\gamma h)(\bb)}{(\aa)^2}\right)\\ 
+\frac{\bar{y}}{N^{1/4}(\mathrm{ch}(\gamma h) - \mathrm{ch}(\beta))}\left(ij +\mathrm{th}(\gamma h)+ik\right)
-\frac{\bar{z}}{N^{1/4}}\left(\frac{ij\mathrm{sh}(\gamma h)\mathrm{th}(\beta)+ik\mathrm{sh}(\beta)}{(\mathrm{ch}(\gamma h) - \mathrm{ch}(\beta))} +  jk\right)\\ + ij\frac{\bar{u}}{N^{1/4}}
+ \frac{\bar{v}}{2N^{1/4}(\aa)}\bigg(\frac{i}{\mathrm{sh}(\gamma h)(\bb)}\\
 - \frac{j}{\mathrm{th}(\gamma h) \mathrm{sh}(\beta)(\bb)}
- ijk\frac{\mathrm{ch}(\gamma h)}{\mathrm{ch}(\beta)(\aa)^2}\bigg)
+ ijk\frac{\bar{w}}{N^{1/4}(\aa)^2}  \bigg]\,.
\end{multline} 
}}
\end{allowdisplaybreaks}
\fine

\begin{COR}\label{lemma Gspace scaling}
In the setting of Proposition \ref{lemmaMP} let consider a function $\varphi(r_N(t)$,$\,\bar{x}_N(t)$,$\,\bar{y}_N(t)$, $\,\bar{z}_N(t)$,$\,\bar{u}_N(t)$,$\,\bar{v}_N(t)$,$\,\bar{w}_N(t))=\varphi(\bar{z}_N(t))$, with $\varphi \in \mathcal{C}^1$.
Then the infinitesimal generator $\mathcal{G}_N$ of the Markov process $(r_N(t)$,$\,\bar{x}_N(t)$,$\,\bar{y}_N(t)$,$\,\bar{z}_N(t)$,$\,\bar{u}_N(t)$,$\,\bar{v}_N(t)$,$\,\bar{w}_N(t))$ applied to that particular function satisfies:
\begin{equation}\label{Gspace scaling}%IG4
\begin{split}
\mathcal{G}_N &\varphi(\bar{z})= \\
&=2\varphi_{\bar{z}}\bigg\{\!\!-\cosh(\gamma h)\,\bar{z}+\frac{\g\cosh(\gamma h)}{N^{1/4}(\cosh(\beta)+\cosh(\gamma h))}\,r\left[\bar{x}+\frac{\bar{v}}{2\sinh(\gamma h)\left(\cosh(\beta)+2\cosh(\gamma h)\right)}\right]\\
&\qquad\qquad-\frac{\gamma}{N^{1/4}\sinh(\beta)\left(\cosh(\beta)+\cosh(\gamma h)\right)^2}\,\bigg[\cosh(\beta)\sinh(\gamma h)\bar{x}^{\,2}\\
&\qquad\qquad\qquad\qquad-\frac{\bar{v}^2}{4\tanh(\gamma h)\left(\cosh(\beta)+2\cosh(\gamma h)\right)^2} +\frac{(\cosh(\beta)-\cosh(\gamma h))}{2(\cosh(\beta)+2\cosh(\gamma h))}\,\bar{x}\bar{v}\bigg]\bigg\}\\
&\qquad+o\left( \frac{1}{N^{1/4}}\right)\,,
\end{split}
\end{equation}
where the remainders are continuous functions of $(r,\bar{x},\bar{y},\bar{z},\bar{u},\bar{v},\bar{w})$ and they are of order $o(1/N^{1/4})$ pointwise, but not uniformly in $(r,\bar{x},\bar{y},\bar{z},\bar{u},\bar{v},\bar{w})$.
\end{COR}

\Proof 
Let consider the function $\varphi(r_N(t),\bar{x}_N(t),\bar{y}_N(t),\bar{z}_N(t),\bar{u}_N(t),\bar{v}_N(t),\bar{w}_N(t))=$ \\ $\varphi(\bar{z}_N(t))$, then \eqref{IG5'} becomes
\begin{multline*}
\mathcal{G}_N \varphi(\bar{z}) =  \sum_{i, j, k\in\mathscr{S}} \vert A_{N}(i, j, k)\vert  e^{\gjm}\,\cdot\\
\quad\cdot 
%&\,\cdot 
\left[\varphi\left(\bar{z} - jk \frac{2}{N^{3/4}}\right)
-\varphi(\bar{z}) \right] \,,
\end{multline*}
where $\vert A_{N}(i, j, k)\vert$ is given by \eqref{contatore}.\\
Now we develop $\varphi$ around $\bar{z}$ with the Taylor expansion stopped at the second order. So,
\begin{multline*}
=\sum_{i, j, k\in\mathscr{S}} \vert A_{N}(i, j, k)\vert  e^{-\gamma jkh}e^{-\g j\left( \frac{\bar{x}}{N^{1/4}(\mathrm{ch}(\beta)+\mathrm{ch}(\gamma h))} + \frac{\bar{v}}{2N^{1/4} \mathrm{sh}(\gamma h)(\mathrm{ch}(\beta)+\mathrm{ch}(\gamma h))(\mathrm{ch}(\beta)+2\mathrm{ch}(\gamma h))}\right)}\,\cdot\\
\quad\cdot 
%&\,\cdot 
\left[-\varphi_{\bar{z}}\frac{2}{N^{3/4}}\,jk+\varphi_{\bar{z}\bar{z}}\frac{2}{N^{3/2}}+o\left(\frac{1}{N^{3/2}}\right)\right]
\end{multline*}
and if we replace the exponential functions with the expression $e^{-\alpha}=\mathrm{ch}(\alpha)-\mathrm{sh}(\alpha)$, by simple computations we obtain that
\begin{allowdisplaybreaks}
\begin{multline}\label{primo sviluppo}
=2\varphi_{\bar{z}}N^{1/4}\bigg\{\mathrm{sh}(\g h)\mathrm{ch}\bigg[\frac{\g}{N^{1/4}(\aa)} \left(\bar{x} + \frac{\bar{v}}{2 \mathrm{sh}(\gamma h)(\bb)}\right)\bigg]\\
+\mathrm{ch}(\g h)\frac{r}{N^{1/2}}\mathrm{sh}\bigg[\frac{\g}{N^{1/4}(\aa)} \left(\bar{x} + \frac{\bar{v}}{2 \mathrm{sh}(\gamma h)(\bb)}\right)\bigg]\\
-\mathrm{sh}(\g h)\left[\frac{1}{N^{1/4}\left(\aa\right)}\left(\frac{ \bar{x}}{\mathrm{th}(\b)}- \frac{\bar{v}}{2\mathrm{sh}(\b)\mathrm{th}(\g h)(\bb)}\right)\right]\cdot\\
\cdot\mathrm{sh}\bigg[\frac{\g}{N^{1/4}(\aa)} \left(\bar{x} + \frac{\bar{v}}{2 \mathrm{sh}(\gamma h)(\bb)}\right)\bigg]\\
-\mathrm{ch}(\g h)\left[\frac{\bar{z}}{N^{1/4}}+\mathrm{th}(\g h)\right]\mathrm{ch}\bigg[\frac{\g}{N^{1/4}(\aa)} \left(\bar{x} + \frac{\bar{v}}{2 \mathrm{sh}(\gamma h)(\bb)}\right)\bigg]\bigg\}\\
+\frac{2\varphi_{\bar{z}\bar{z}}}{N^{1/2}}\bigg\{\mathrm{ch}(\g h)\mathrm{ch}\left[\frac{\g}{N^{1/4}(\aa)} \left(\bar{x} + \frac{\bar{v}}{2 \mathrm{sh}(\gamma h)(\bb)}\right)\right]\\
+\mathrm{sh}(\g h)\frac{r}{N^{1/2}}\mathrm{sh}\left[\frac{\g}{N^{1/4}(\aa)} \left(\bar{x} + \frac{\bar{v}}{2 \mathrm{sh}(\gamma h)(\bb)}\right)\right]\\
-\mathrm{ch}(\g h)\left[\frac{1}{N^{1/4}\left(\aa\right)}\left(\frac{ \bar{x}}{\mathrm{th}(\b)}- \frac{\bar{v}}{2\mathrm{sh}(\b)\mathrm{th}(\g h)(\bb)}\right)\right]\cdot\\
\cdot\mathrm{sh}\bigg[\frac{\g}{N^{1/4}(\aa)} \left(\bar{x} + \frac{\bar{v}}{2 \mathrm{sh}(\gamma h)(\bb)}\right)\bigg]\\
-\mathrm{sh}(\g h)\left[\frac{\bar{z}}{N^{1/4}}+\mathrm{th}(\g h)\right]\mathrm{ch}\bigg[\frac{\g}{N^{1/4}(\aa)} \left(\bar{x} + \frac{\bar{v}}{2 \mathrm{sh}(\gamma h)(\bb)}\right)\bigg]\bigg\}=(\ast)\,.
\end{multline}
\end{allowdisplaybreaks}\\
In what follows we will also consider the Taylor expansions stopped at the second order of the following terms
\begin{equation}\label{sinh}
\begin{split}
\mathrm{sh} &\left[\frac{\g}{N^{1/4}\left(\aa\right)}\left(\bar{x}+\frac{\bar{v}}{2\mathrm{sh}(\g h)\left(\bb\right)}\right)\right] =\\
&=\frac{\g}{N^{1/4}\left(\aa\right)}\left(\bar{x}+\frac{\bar{v}}{2\mathrm{sh}(\g h)\left(\bb\right)}\right)+o\left(\frac{1}{N^{1/2}}\right)
\end{split}
\end{equation}
\begin{equation}\label{cosh}
\begin{split}
\mathrm{ch} &\left [\frac{\g}{N^{1/4}\left(\aa\right)}\left(\bar{x}+\frac{\bar{v}}{2\mathrm{sh}(\g h)\left(\bb\right)}\right)\right]=\\
&=1+\frac{\g^2}{2N^{1/2}\left(\aa\right)^2}\left(\bar{x}+\frac{\bar{v}}{2\mathrm{sh}(\g h)\left(\bb\right)}\right)^2+ o\left(\frac{1}{N^{1/2}}\right)\,.
\end{split}
\end{equation}
Then, if we consider the previous expansions and we reorder the terms, we find that 
\begin{allowdisplaybreaks}
\begin{multline*}
(\ast)=
2\varphi_{\bar{z}}\bigg\{\!\!-\mathrm{ch}(\g h)\,\bar{z}+\frac{\g\mathrm{ch}(\g h)}{N^{1/4}(\aa)}\,r\left[\bar{x}+\frac{\bar{v}}{2\mathrm{sh}(\g h)(\bb)}\right]\\
-\frac{\g}{N^{1/4}\mathrm{sh}(\b)(\aa)^2}\,\bigg[\mathrm{ch}(\b)\mathrm{sh}(\g h)\bar{x}^{\,2}\\
-\frac{\bar{v}^2}{4\mathrm{th}(\g h)(\bb)^2}+\frac{(\mathrm{ch}(\b)-\mathrm{ch}(\g h))}{2(\bb)}\,\bar{x}\bar{v}\bigg]\bigg\}
+o\left( \frac{1}{N^{1/4}}\right)\,,
\end{multline*}
\end{allowdisplaybreaks}\\
which is just \eqref{Gspace scaling}.
\fine

We define
\begin{multline*}
(r_N(t), \widetilde{x}_N(t), \widetilde{y}_N(t)\, \widetilde{z}_N(t), \widetilde{u}_N(t), \widetilde{v}_N(t), \widetilde{w}_N(t)) := \\
:= (r_{N}(t), \bar{y}_{N}(N^{1/4}t), \bar{z}_{N}(N^{1/4}t), \bar{u}_{N}(N^{1/4}t), \bar{v}_{N}(N^{1/4}t), \bar{w}_{N}(N^{1/4}t))
\end{multline*}
and we consider the infinitesimal generator, $\mathcal{J}_N=N^{1/4}\mathcal{G}_N$, subject to the time-rescaling and applied to the particular function $$\varphi(r_N(t), \widetilde{x}_N(t), \widetilde{y}_N(t)\, \widetilde{z}_N(t), \widetilde{u}_N(t), \widetilde{v}_N(t), \widetilde{w}_N(t))=(\widetilde{z}_N(t))^2\,.$$ 
We choose this kind of function, since $(\widetilde{z}_N(t))^2$ is a sequence of positive semimartingales on a suitable probability space $(\Omega, \mathcal{A},P)$ and then the following decomposition holds:
\begin{equation}\label{decompz} (\widetilde{z}_N(t))^2-(\widetilde{z}_N(0))^2=\int_0^t\mathcal{J}_N(\widetilde{z}_N(s))^2\,ds+\mathcal{M}_{N,\widetilde{z}^2}^t\,. 
\end{equation}
% \[ d(\widetilde{z}_N(t))^2=\mathcal{J}_N(\widetilde{z}_N(t))^2\,dt+d\mathcal{M}_{N,\widetilde{z}^2}^t  \]
In (\ref{decompz}):
\begin{equation*}
\begin{split}
\mathcal{J}_N(\widetilde{z}_N(t))^2=& N^{1/4} 
\sum_{i,j,k\in \mathscr{S}}\left|A(i,j,k,N^{1/4}t)\right|\,\,\cdot\\
&\,\cdot
e^{\gjm}\overline{\nabla}^{(j)}\left[(\widetilde{z}_N(s))^2\right]\,,
\end{split}
\end{equation*}
\[\mathcal{M}_{N,\widetilde{z}^2}^t=\int_0^t\sum_{i,j,k\in \mathscr{S}}\overline{\nabla}^{(j)}\left[(\widetilde{z}_N(s))^2\right]
\,\widetilde{\Lambda}^\omega_N(i,j,k,ds)\,, \]
which is a local martingale, where 
\begin{equation}\label{gradsemimart}
\overline{\nabla}^{(j)}\left[(\widetilde{z}_N(s))^2\right]:=
\left[\left(\widetilde{z}_N(s)- jk\frac{2}{N^{3/4}} \right)^2- (\widetilde{z}_N(s))^2\right]\,
\end{equation}
and
\begin{equation}\label{lambdatilde}
\widetilde{\Lambda}^\omega_N(i,j,k,dt) := \Lambda^\omega_N(i,j,k,dt)-\underbrace{N^{1/4}\left\vert A(i,j,k,N^{1/4}t)\right\vert  e^{-\gamma j\left(m^\sigma_{\rho_N(N^{1/4}t)}+kh\right)} dt}_{:= \lambda(i,j,k,t)\,dt}\,.
\end{equation}
The counter $\left\vert A(i,j,k,N^{1/4}t)\right\vert$ is given in analogy with \eqref{contatore}, replacing the variables $r$, $\bar{x}$, $\bar{y}$, $\bar{z}$, $\bar{u}$, $\bar{v}$, $\bar{w}$ with the stochastic processes $r_N(t)$, $\widetilde{x}_N(t)$, $\widetilde{y}_N(t)$, $\widetilde{z}_N(t)$, $\widetilde{u}_N(t)$, $\widetilde{v}_N(t)$ and $\widetilde{w}_N(t)$.\\
%by
%\begin{multline*}
%\vert A (i, j, k, N^{1/4}t) \vert = \frac{N}{8} \bigg[1 + k m^{\underline{\eta}}_{\rho_N(N^{1/4}t)} + i m^{\underline{\sigma}}_{\rho_N(N^{1/4}t)} + j m^{\underline{\omega}}_{\rho_N(N^{1/4}t)} + ij m^{\underline{\sigma} \, \underline{\omega}}_{\rho_N(N^{1/4}t)} \\
%+ ik m^{\underline{\sigma} \, \underline{\eta}}_{\rho_N(N^{1/4}t)} + jk m^{\underline{\omega} \, \underline{\eta}}_{\rho_N(N^{1/4}t)} + ijk m^{\underline{\sigma} \, \underline{\omega} \, \underline{\eta}}_{\rho_N(N^{1/4}t)}\bigg] \,.
%\end{multline*}
%%%%%%%%%%%%%%%%%%%%%%%%%%%%%%%%%%%%%%%%%%%%%%%%%%%%%%%%%%%%%%%%%%%%%%%%%%%%%%%%%%%%%%%%%%%%%%%%%
% \[ d(\widetilde{y}_N(t))^2=\mathcal{J}_N(\widetilde{y}_N(t))^2\,dt+d\mathcal{M}_{N,\widetilde{y}^2}^t  \]
As we can evidently see, $\widetilde{\Lambda}^\omega_N(i,j,k,dt)$ is the difference between the point process
 $\Lambda^\omega_N(i,j,k,dt)$, defined on \mbox{$\mathscr{S}^3\times\mathbb{R}^+$}, and its intensity $\lambda(i,j,k,t)\,dt$. 

\begin{REM}
If we call $(\mathcal{A}_t)_{t\geq0}$ a filtration generated by $\Lambda^\omega_N$ on $(\Omega,\mathcal{A},P)$, then the processes $\mathcal{J}_N(\widetilde{z}_N(t))^2$ and  $\overline{\nabla}^{(j)}[(\widetilde{z}_N(t))^2]$ are $\mathcal{A}_t-$adapted processes.
\end{REM} 

As a consequence of the considerations just explained, we are in the proper situation to use a result about collapsing processes appeared in \cite{CoEi88} and then slightly generalized in \cite{Sar07}. Here we recall it in Proposition \ref{collapse} and then we adapt it to our specific case in Lemma \ref{LmmClps}, that is what we have to prove.

\begin{PROP}\label{collapse}
\textsf{Let $\{X_n(t)\}_{n \geq 1}$ be a sequence of positive semimartingales on a probability space $(\Omega, \mathscr{A}, \mathscr{P})$, with
\[dX_n(t) = S_n(t)dt + \int_{\mathbb{R}^+\times \mathscr{Y}} f_n(t^-, y) [\Lambda_n(dt,dy) - A_n(t,dy)dt]\,.\]
Here, $\Lambda_n$ is a point-process of intensity $A_n(t,dy)dt$ on $\mathbb{R}^+\times \mathscr{Y}$, where $\mathscr{Y}$ is a measurable space, and $S_n(t)$ and $f_n(t)$ are $\mathscr{A}_t$-adapted processes, if we consider $(\mathscr{A}_t)_{t\geq 0}$ a filtration on $(\Omega, \mathscr{A}, \mathscr{P})$ generated by $\Lambda_n$.\\ 
Let $d>1$ and $C_i$ constants independent of $n$ and $t$. Suppose there exist $\{\alpha_n\}_{n \geq1}$ and $\{\beta_n\}_{n \geq 1}$, increasing sequences with
\[n^{1/d}\alpha_n^{-1} \xrightarrow{n\rightarrow +\infty} 0, n^{-1} \alpha_n  \xrightarrow{n\rightarrow +\infty} 0, n^{-1}\beta_n  \xrightarrow{n\rightarrow +\infty} 0  \]
and
\[ E \Big[ \Big( X_n(0) \Big)^d \Big] \leq C_1 \alpha_n^{-d} \qquad \mbox{for all $n$}\,. \]
Furthermore, let $\{\tau_n\}_{n \geq 1}$ be stopping times such that, for $t \in [0, \tau_n]$ and $n \geq 1$,
\[ S_n(t) \leq -n \delta X_n(t) + \beta_n C_2 + C_3 \qquad \mbox{with $\delta>0$}\,, \] 
\[ \sup_{\omega\in\Omega, y\in\mathscr{Y}, t\leq\tau_n}\vert f_n(t,y) \vert \leq C_4 \alpha_{n}^{-1}\,.\]
Hence:
\begin{enumerate}
\item if it holds
\begin{equation}\label{condthmCOEI}\tag{$\star\star$}
\int_{\mathscr{Y}} (f_n(t,y))^2 A_n(t,dy) \leq C_5\,,
\end{equation}
then, for any $\varepsilon >0$, there exist $C_6>0$ and $n_0$ such that
\begin{equation}\label{tesithmCOEI}\tag{$\star\star\star$}
\sup_{n \geq n_0} \mathscr{P} \left\{ \sup_{0 \leq t \leq T \wedge \tau_n} X_n (t) > C_6 (n^{1/d} \alpha_n^{-1} \vee \alpha_n n^{-1})\right\} \leq \varepsilon\,; 
\end{equation}
\item if instead of \eqref{condthmCOEI} we have
\[ \int_{\mathscr{Y}} (f_n(t,y))^2 A_n(t,dy) \leq C_5(X_n(t) + n^{-1})\,,\]
then, instead of \eqref{tesithmCOEI}, we get
\[\sup_{n \geq n_0} \mathscr{P} \left\{ \sup_{0 \leq t \leq T \wedge \tau_n} X_n (t) > C_6 (n^{1/d} \alpha_n^{-1} \vee \beta_n n^{-1})\right\} \leq \varepsilon\,. \]
\end{enumerate}}
\end{PROP}

\begin{LEMMA}\label{LmmClps}
Consider $d>2$, $\delta > 0$  and $\kappa := \kappa(N)$, such that $\kappa \xrightarrow{N\rightarrow +\infty} +\infty$. For $t \in [0, \tau_N^M]$ and $N \geq 1$, there exist constants $C_{\cdot}$'s independent of $N$ and $t$ and two increasing sequences  $\{\alpha_N\}_{N \geq 1}$ and $\{\beta_N\}_{N \geq 1}$, which satisfy the following conditions:
\begin{equation}\label{C1}
   \kappa^{1/d}\alpha_N^{-1}\xrightarrow{N\rightarrow +\infty} 0,\kappa^{-1}\alpha_N \xrightarrow{N\rightarrow +\infty} 0,\kappa^{-1}\beta_N \xrightarrow{N\rightarrow +\infty} 0 \, ,
\end{equation}
\begin{equation}\label{C2}
    E\left[(\widetilde{z}_N(0))^{2d}\right]\leq C_{1}\,\alpha_N^{-2d} \quad \mbox{for all $N$,}
\end{equation}
\begin{equation}\label{C3}
    \mathcal{J}_N(\widetilde{z}_N(t))^2 \leq -\kappa\delta(\widetilde{z}_N(t))^2+\beta_N C_{2}+C_{3} \, ,
\end{equation}
\begin{equation}\label{C4}
    \sup_{\omega\in\Omega, j,k\in\mathscr{S}, t\leq\tau_N} \left\vert \overline{\nabla}^{(j)}\left[(\widetilde{z}_N(t))^2\right]  \right\vert \leq C_{4}\, \alpha_N^{-1} \, ,
\end{equation}
\begin{equation}\label{C5}
    \sum_{i,j,k\in\mathscr{S}}\Big[ \overline{\nabla}^{(j)}\left[(\widetilde{z}_N(t))^2\right] \Big]^2\,\lambda(i,j,k,t) \leq C_{5}\left((\widetilde{z}_N(t))^2   + \kappa^{-1}\right) \, 
\end{equation}
and such that, for every $\varepsilon>0$, the following estimate holds
\begin{equation}\label{thLmmClps}
    \sup_{N\geq N_{0}} P\left\{\sup_{0\leq t\leq T\wedge\tau_N^M}(\widetilde{z}_N(t))^2\stackrel{(*)}{>}C_{6}\,\left(\kappa^{1/d}\alpha_N^{-1} \vee \kappa^{-1}\beta_N \right)\right\}\leq\varepsilon\,.
\end{equation}
\end{LEMMA}

\Proof
We aim to prove the existence of these sequences $\{\alpha_N\}_{N \geq 1}$, $\{\beta_N\}_{N \geq 1}$ and constants $C_{\cdot}$ and to give a characterization of them. We show that the hypotheses required by Lemma \ref{LmmClps} hold true.\\

\eqref{C2}:  From  \eqref{fluctcrt} we get
\[ 
\widetilde{z}_N(0) = N^{1/4} \left(m^{\underline{\o} \, \underline{\eta}}_{\rho_N(0)}-\tanh(\g h)\right) \,.\]
The random variables $(\o_j(0),\eta_j)$ are independent, so a Central Limit Theorem applies: in the limit as $N \longrightarrow +\infty$, 
\[
N^{1/4}\widetilde{z}_N(0) =  N^{1/2}\left(m^{\underline{\o} \, \underline{\eta}}_{\rho_N(0)}-\tanh(\g h)\right)
\] 
converges to a Gaussian random variable and, since $m^{\underline{\o} \, \underline{\eta}}_{\rho_N(0)} \in [-1,+1]$, there is convergence of all the moments. Thus,
\[
E\left[N^d\left(m^{\underline{\o} \, \underline{\eta}}_{\rho_N(0)}-\tanh(\g h)\right)^{2d}\right] \leq C_{1}
\]
and we obtain the following estimate for the $2d$-th moments of $\widetilde{z}_N(0)$: 
\begin{eqnarray*}
E\left[ \left( \widetilde{z}_N(0)\right)^{2d} \right] &=&  E\left[ N^{d/2} \left(m^{\underline{\o} \, \underline{\eta}}_{\rho_N(0)}-\tanh(\g h)\right)^{2d} \right]\\
&=& N^{-d/2} \, E\left[ N^{d} \left(m^{\underline{\o} \, \underline{\eta}}_{\rho_N(0)}-\tanh(\g h)\right)^{2d} \right]  \leq   C_{1} \,  N^{-d/4}\,.
\end{eqnarray*}
Thus \eqref{C2} holds.\\

\eqref{C3}: For $t\in [0,\tau_N^M]$ we consider the Taylor expansions of the hyperbolic sine and cosine stopped at the second order (see \eqref{sinh} and \eqref{cosh}) with the Lagrangian expressions of their remainders estimated as follows: 
\begin{allowdisplaybreaks}
\begin{align}\label{resto sinh}
\vert R_s\vert & \leq\sup\bigg\{\vert\mathrm{ch}(\theta)\vert:
\theta\in \left[0,\frac{\g}{N^{1/4}\left(\aa\right)}\left(\widetilde{x}_N(t)+\frac{\widetilde{v}_N(t)}{2\mathrm{sh}(\g h)\left(\bb\right)}\right)\right]\bigg\}\cdot\nonumber\\
&\qquad\cdot\frac{\g^3}{6N^{3/4}\left(\aa\right)^3}\left(\widetilde{x}_N(t)+\frac{\widetilde{v}_N(t)}{2\mathrm{sh}(\g h)\left(\bb\right)}\right)^3\nonumber\\ 
& \leq\mathrm{ch}\left[\frac{\g M}{N^{1/4}\left(\aa\right)}\left(1+\frac{1}{2\mathrm{sh}(\g h)\left(\bb\right)}\right)\right]\cdot\nonumber\\
&\qquad\cdot\frac{\g^3 M^3}{6N^{3/4}(\aa)^3}\left(1+\frac{1}{2\mathrm{sh}(\g h)(\bb)}\right)^3
\end{align}
\begin{align}\label{resto cosh}
\vert R_c\vert & \leq\sup\bigg\{\vert\mathrm{sh}(\theta)\vert:
\theta\in \left[0,\frac{\g}{N^{1/4}\left(\aa\right)}\left(\widetilde{x}_N(t)+\frac{\widetilde{v}_N(t)}{2\mathrm{sh}(\g h)\left(\bb\right)}\right)\right]\bigg\}\cdot\nonumber\\
&\qquad\cdot\frac{\g^3}{6N^{3/4}\left(\aa\right)^3}\left(\widetilde{x}_N(t)+\frac{\widetilde{v}_N(t)}{2\mathrm{sh}(\g h)\left(\bb\right)}\right)^3\nonumber\\ 
& \leq\mathrm{sh}\left[\frac{\g M}{N^{1/4}\left(\aa\right)}\left(1+\frac{1}{2\mathrm{sh}(\g h)\left(\bb\right)}\right)\right]\cdot\nonumber\\
&\qquad\cdot\frac{\g^3 M^3}{6N^{3/4}(\aa)^3}\left(1+\frac{1}{2\mathrm{sh}(\g h)(\bb)}\right)^3\,.
\end{align}
\end{allowdisplaybreaks}

Now we derive the particular characterization of $\mathcal{J}_N(\widetilde{z}_N(t))^2$, adapting the explicit expression of $\mathcal{G}_N\varphi\left(\widetilde{z}_N(t)\right)$ found in Corollary \ref{lemma Gspace scaling} (in other words, setting $\varphi(\widetilde{z}_N(t))=(\widetilde{z}_N(t))^2$, and taking into account the time-rescaling). Then, we proceed to find an upper bound for this quantity. Thus, by \eqref{primo sviluppo}, \eqref{sinh} and \eqref{cosh}, if we reorder the terms, we get
\begin{allowdisplaybreaks}
\begin{multline*}
\mathcal{J}_N
(\widetilde{z}_N(t))^2 =%4N^{1/2} \widetilde{z}_N(t)\cdot\\ 
-4N^{1/4}\mathrm{ch}(\g h)\left(\widetilde{z}_N(t)\right)^2
-\frac{4\g \mathrm{sh}(\g h)}{\mathrm{th}(\b)(\aa)^2}\left(\widetilde{x}_N(t)\right)^2\widetilde{z}_N(t)\\
-\frac{2\g(\mathrm{ch}(\b)-\mathrm{ch}(\g h))}{\mathrm{sh}(\b)(\aa)^2(\bb)}\,\widetilde{x}_N(t)\,\widetilde{v}_N(t)\,\widetilde{z}_N(t)\\
+\frac{4\g}{\mathrm{sh}(\b)\mathrm{th}(\g h)(\aa)^2(\bb)^2}\left(\widetilde{v}_N(t)\right)^2\widetilde{z}_N(t)\\
+\frac{4\g\mathrm{ch}(\g h)}{\aa}\,\frac{r_N(t)\,\widetilde{z}_N(t)}{N^{1/4}}\left[\widetilde{x}_N(t)+\frac{\widetilde{v}_N(t)}{\mathrm{sh}(\g h)\left(\bb\right)}\right]\\
-\mathrm{ch}(\g h)\,\frac{\left(\widetilde{z}_N(t)\right)^2}{N^{1/4}}\bigg[\frac{2\g^2}{(\aa)^2}\left(\widetilde{x}_N(t)\right)^2\\
+\frac{2\g^2}{\mathrm{sh}(\g h)(\aa)^2(\bb)}\,\widetilde{x}_N(t)\,\widetilde{v}_N(t)\\
+\frac{\g^2}{2\mathrm{sh}^2(\g h)(\aa)^2(\bb)^2}\left(\widetilde{v}_N(t)\right)^2\bigg]\\
+R_s\bigg[4\mathrm{ch}(\g h)\,r_N(t)\widetilde{z}_N(t)-4N^{1/4}\frac{\mathrm{sh}(\g h)}{\mathrm{th}(\b)(\aa)}\,\widetilde{x}_N(t)\,\widetilde{z}_N(t)\\
+2N^{1/4}\frac{\mathrm{ch}(\g h)}{\mathrm{sh}(\b)(\aa)(\bb)}\,\widetilde{v}_N(t)\,\widetilde{z}_N(t)\bigg]\\
+R_c\bigg[-4N^{1/4}\mathrm{ch}(\g h)\left(\widetilde{z}_N(t)\right)^2+\frac{4\mathrm{ch}(\g h)}{N^{1/4}}\bigg]\\
%%%%%%%%%%%%%%%%%%
\\
\leq -4N^{1/4}\mathrm{ch}(\g h)\left(\widetilde{z}_N(t)\right)^2
+\frac{4\g \mathrm{sh}(\g h)}{\mathrm{th}(\b)(\aa)^2}\left(\widetilde{x}_N(t)\right)^2\vert\widetilde{z}_N(t)\vert\\
+\frac{2\g}{\mathrm{sh}(\b)(\aa)(\bb)}\,\vert\widetilde{x}_N(t)\vert\,\vert\widetilde{v}_N(t)\vert\,\vert\widetilde{z}_N(t)\vert\\
+\frac{4\g}{\mathrm{sh}(\b)\mathrm{th}(\g h)(\aa)^2(\bb)^2}\left(\widetilde{v}_N(t)\right)^2\vert\widetilde{z}_N(t)\vert\\
+\frac{4\g\mathrm{ch}(\g h)}{\aa}\,\vert r_N(t)\vert\,\vert\widetilde{z}_N(t)\vert\left[\vert\widetilde{x}_N(t)\vert+\frac{\vert\widetilde{v}_N(t)\vert}{\mathrm{sh}(\g h)\left(\bb\right)}\right]\\
+\mathrm{ch}(\g h)\,\left(\widetilde{z}_N(t)\right)^2\bigg[\frac{2\g^2}{(\aa)^2}\left(\widetilde{x}_N(t)\right)^2\\
+\frac{2\g^2}{\mathrm{sh}(\g h)(\aa)^2(\bb)}\,\vert\widetilde{x}_N(t)\vert\,\vert\widetilde{v}_N(t)\vert\\
+\frac{\g^2}{2\mathrm{sh}^2(\g h)(\aa)^2(\bb)^2}\left(\widetilde{v}_N(t)\right)^2\bigg]\\
+\vert R_s\vert\bigg[4\mathrm{ch}(\g h)\,\vert r_N(t)\vert \,\vert\widetilde{z}_N(t)\vert+4N^{1/4}\frac{\mathrm{sh}(\g h)}{\mathrm{th}(\b)(\aa)}\,\vert\widetilde{x}_N(t)\vert\,\vert\widetilde{z}_N(t)\vert\\
+2N^{1/4}\frac{\mathrm{ch}(\g h)}{\mathrm{sh}(\b)(\aa)(\bb)}\,\vert\widetilde{v}_N(t)\vert\,\vert\widetilde{z}_N(t)\vert\bigg]\\
+\vert R_c\vert\bigg[4N^{1/4}\mathrm{ch}(\g h)\left(\widetilde{z}_N(t)\right)^2+4\mathrm{ch}(\g h)\bigg]\\
%%%%%%%%%%%%%%%%%%
\\
\leq
-4N^{1/4}\mathrm{ch}(\g h)\left(\widetilde{z}_N(t)\right)^2
+4M^3\frac{\g \mathrm{sh}(\g h)}{\mathrm{th}(\b)(\aa)^2}\\
+2M^3\frac{\g}{\mathrm{sh}(\b)(\aa)(\bb)}\\
+4M^3\frac{\g}{\mathrm{sh}(\b)\mathrm{th}(\g h)(\aa)^2(\bb)^2}\\
+4M^3\frac{\g\mathrm{ch}(\g h)}{\aa}\left[1+\frac{1}{2\mathrm{sh}(\g h)\left(\bb\right)}\right]\\
+M^4\mathrm{ch}(\g h)\bigg[\frac{2\g^2}{(\aa)^2}
+\frac{2\g^2}{\mathrm{sh}(\g h)(\aa)^2(\bb)}\\
+\frac{\g^2}{2\mathrm{sh}^2(\g h)(\aa)^2(\bb)^2}\bigg]\\
%%%%%%%%%%%%resti
+\mathrm{ch}\left[\frac{\g M}{\left(\aa\right)}\left(1+\frac{1}{2\mathrm{sh}(\g h)\left(\bb\right)}\right)\right]\frac{\g^3 M^3}{6(\aa)^3}\cdot\\
\cdot\left(1+\frac{1}{2\mathrm{sh}(\g h)(\bb)}\right)^3\bigg[4M^2\mathrm{ch}(\g h)+4M^2\frac{\mathrm{sh}(\g h)}{\mathrm{th}(\b)(\aa)}\\
+2M^2\frac{\mathrm{ch}(\g h)}{\mathrm{sh}(\b)(\aa)(\bb)}\bigg]\\
+\mathrm{sh}\left[\frac{\g M}{\left(\aa\right)}\left(1+\frac{1}{2\mathrm{sh}(\g h)\left(\bb\right)}\right)\right]\frac{\g^3 M^3}{6(\aa)^3}\cdot\\
\cdot\left(1+\frac{1}{2\mathrm{sh}(\g h)(\bb)}\right)^3\bigg[4M^2\mathrm{ch}(\g h)+4\mathrm{ch}(\g h)\bigg]\\
\\
= -4N^{1/4}\mathrm{ch}(\g h)(\widetilde{z}_N(t))^2 + C_{2} + C_{3}
\end{multline*}
\end{allowdisplaybreaks}\\
Hence, we have obtained the desired inequality if we choose:  \mbox{$\kappa := N^{1/4}$}, \mbox{$\delta := 4 \mathrm{ch} (\g h)$} (which is a positive constant), \mbox{$\beta_N \equiv 1$} and \mbox{$C_{2} + C_{3}$} equal to the rest of the expression, which is constant with respect to $N$ and $t$, as required.\\

\eqref{C4}: Now, we evaluate the supremum of the modulus of $\overline{\nabla}^{(j)}\left[(\widetilde{z}_N(s))^2\right]$, defined as in \eqref{gradsemimart}. It easily yields
\begin{align*}
\sup_{\omega\in\Omega, j,k\in\mathscr{S}, t\in[0,\tau_N^M]} \left\vert \overline{\nabla}^{(j)}\left[(\widetilde{z}_N(t))^2\right]  \right\vert&= \sup_{\omega\in\Omega, j,k\in\mathscr{S}, t\in[0,\tau_N^M]} \bigg\vert \frac{4}{N^{3/2}}-jk\frac{4\widetilde{z}_N(t)}{N^{3/4}}\bigg\vert \\
&\leq \frac{4}{N^{5/8}}(1+M)\,N^{-1/8} \leq C_{4}\,N^{-1/8}\,,
\end{align*}
where we set $C_{4}=4(1+M)$ and $\alpha_N=N^{1/8}$.\\

\eqref{C5}: Recalling the definition of $\overline{\nabla}^{(j)}[(\widetilde{z}_N(s))^2]$ and of $\lambda(i,j,k,t)$, which we can find in \eqref{gradsemimart} and in \eqref{lambdatilde}, we have
\begin{allowdisplaybreaks}
\begin{multline*}
N^{1/4}  \sum_{i,j,k\in\mathscr{S}} \vert A(i,j,k,N^{1/4}t)\,\vert \, e^{\gjmm}\,\cdot\\ 
\cdot\bigg[ \bigg(\widetilde{z}_N(t) -  jk\frac{2}{N^{3/4}}\bigg)^2  - (\widetilde{z}_N(t))^2  \bigg]^2 \\
= N^{1/4} \sum_{i,j,k\in\mathscr{S}} \vert A(i,j,k,N^{1/4}t)\,\vert \, e^{\gjmm}\,\cdot\\
\cdot\bigg[\frac{4}{N^{3/2}}  -jk\frac{4}{N^{3/4}}\,\widetilde{z}_N(t)\bigg]^2 \\
= N^{1/4} \sum_{i,j,k\in\mathscr{S}} \vert A(i,j,k,N^{1/4}t)\,\vert\,  e^{-\gamma jkh}e^{- \frac{\g j}{N^{1/4}(\mathrm{ch}(\beta)+\mathrm{ch}(\gamma h))}\left(\widetilde{x}_N(t) + \frac{\widetilde{v}_N(t)}{2 \mathrm{sh}(\gamma h)(\mathrm{ch}(\beta)+2\mathrm{ch}(\gamma h))}\right)}\,\cdot\\
\cdot\bigg[\frac{16}{N^3}-\frac{32}{N^{9/4}} \widetilde{z}_N{(t)}+\frac{16}{N^{3/2}}(\widetilde{z}_N(t))^2\bigg]\,.
\end{multline*}
\end{allowdisplaybreaks}\\
To find an upper bound for this last quantity, we replace the exponential functions with the expression $e^{-\alpha}=\mathrm{ch}(\alpha)-\mathrm{sh}(\alpha)$ and then we proceed in the same way we previously proved \eqref{C3}. We consider the Taylor expansions of the hyperbolic sine and cosine stopped at the second order (see \eqref{sinh}, \eqref{cosh}) and we estimate their remainders as in \eqref{resto sinh}, \eqref{resto cosh}. So, we can show that for $t \in [0,\tau_N^M]$ it holds
\[
\sum_{i,j,k\in\mathscr{S}}\Big[ \overline{\nabla}^{(j)}\left[(\widetilde{z}_N(t))^2\right] \Big]^2\,\lambda(i,j,k,t) \leq C_5 \left[\left(\widetilde{z}_N(t)\right)^2+\frac{1}{N^{1/4}}\right] \,,
\]
where $C_5$ is a positive constant depending on $M$.

\eqref{C1}: It remains to show that the sequences we have found satisfy the conditions about the convergence to zero. But, 
\[ \lim_{N\rightarrow+\infty}(N^{1/4})^{1/d}(N^{1/8})^{-1}=\lim_{N\rightarrow+\infty}N^{1/(4d)-1/8}=0\quad\Longleftrightarrow\quad d>2 \]
\[ \lim_{N\rightarrow+\infty}N^{1/8}N^{-1/4} = \lim_{N\rightarrow+\infty}N^{-1/8} = 0 \]
\[  \lim_{N\rightarrow+\infty}N^{-1/4} = 0 \]
and hence we know that \eqref{thLmmClps} hold true from Proposition \ref{collapse}. 
\fine

\begin{COR}
We consider the same setting as in Lemma \ref{LmmClps}. For every $\varepsilon>0$ there exists constants $C_{7}$ and $N_{0}$ such that
\begin{equation}\label{thcorclps}
    \sup_{N\geq N_{0}} P\left\{\sup_{0\leq t\leq T\wedge\tau_N^M}\vert\widetilde{z}_N(t)\vert>C_{7}\,\left(\kappa^{1/(2d)}\alpha_N^{-1/2} \vee  \kappa^{-1/2}\beta_N^{1/2}\right)\right\}\leq\varepsilon\,.
\end{equation}
\end{COR}

\Proof 
We set $C_{7}=(C_{6})^{1/2}$ and we extract the square root of the inequality $(*)$ in the previous Lemma to obtain an equivalent set, described in \eqref{thcorclps}, for which the same property holds.
\fine

\begin{REM}
Notice that if we insert the quantities we choose during the proof of Lemma \ref{LmmClps} into \eqref{thcorclps},  we have shown that the following inequality holds
\begin{equation}\label{thcorclps2}
 \sup_{N\geq N_{0}} P\left\{\sup_{0\leq t\leq T\wedge\tau_N^M}\vert\widetilde{z}_N(t)\vert>C_{7}\,\left(N^{1/(8d) - 1/16} \vee  N^{-1/8}\right)\right\}\leq\varepsilon\,.
\end{equation}
\end{REM}

%%%%%%%%%%%%%%%%%%%%%%STEP 3%%%%%%%%%%%%%%%%%%%%%%%%%%
\textbf{STEP 3.} In order to conclude the first part of the proof of Theorem \ref{thmCRTDYN'} we have to show that, for every $\varepsilon >0$ and $N \geq 1$, there exists a constant $M>0$ such that  
\[ P\left\{\tau_N^M\leq T \right\} \leq \epsilon \,.\]
This fact implies that the processes $\widetilde{y}_{N}(t), \widetilde{z}_{N}(t), \widetilde{u}_{N}(t), \widetilde{v}_{N}(t), \widetilde{w}_{N}(t)$ converge to zero in probability, as $N\rightarrow +\infty$, for all $t \in [0,T]$. But, before proving this fact we need the following technical Lemma. 
%%%%%%%%%%%%%%%%%%%%%%%%%%%%%
%%%%%%%%%%%%%%%%%%%%%%%%%%
\begin{LEMMA}\label{LmmSPSC'}
For $t \in [0,T \wedge \tau_{N}^M]$, if we consider $\psi \in \mathcal{C}^1$, a function of the pair of processes
\begin{equation}\label{crtdir'}
\begin{array}{ccl}
r_N(t) & = & N^{1/2}\,  m^{\underline{\eta}}_{\rho_N(t)}\\
\bar{x}_N(t) & = & N^{1/4} \left[\cosh(\gamma h) m^{\underline{\sigma}}_{\rho_N(t)} + \sinh(\beta) m^{\underline{\omega}}_{\rho_N(t)} \right]
\end{array}
\end{equation}
only rescaled in space, then \eqref{IG5'} reduces to
\begin{align}\label{IG4'}
\mathcal{G}_N \psi(r, \bar{x}) =& 2 \bigg[ r \, \frac{\sinh(\beta) \sinh(\gamma h)}{N^{1/4}} - \frac{\gamma^2 \bar{x}^3}{2} \frac{\cosh(\beta) \cosh(\gamma h)}{N^{1/2} (\cosh(\beta) + \cosh(\gamma h))^3}\nonumber \\
                                 & \qquad + \frac{\gamma^3 \bar{x}^3}{6} \frac{\sinh(\beta)(\cosh(\gamma h) - \sinh(\gamma h) \tanh(\gamma h))}{N^{1/2} (\cosh(\beta) + \cosh(\gamma h))^3}\\
                                 & \qquad + \frac{\gamma^2 r \bar{x}^2}{2} \frac{\sinh(\beta) \sinh(\gamma h)}{N^{3/4} (\cosh(\beta) + \cosh(\gamma h))^2} \bigg] \psi_{\bar{x}} + o\bigg( \frac{1}{N^{1/4}}\bigg)\,, \nonumber
\end{align}
where the remainder is a continuous function of $\bar{x}$ and it is of order $o(1/N^{1/4})$ pointwise, but not uniformly in $\bar{x}$.
\end{LEMMA}

\Proof
By \eqref{IG5'}, considering a function $\psi: \mathbb{R}^2 \longrightarrow \mathbb{R}$, $\psi \in \mathcal{C}_b^3$, we deduce
\begin{allowdisplaybreaks}
\begin{align*}
\mathcal{G}_{N}\psi(r, \bar{x}) = & \sum_{i, j, k \in \mathscr{S}} \vert A_{N}(i,j,k) \vert \bigg\{ e^{-\beta ij} \bigg[ \psi\bigg(r, \bar{x} - i \frac{2}{N^{3/4}} \mathrm{ch}(\gamma h) \bigg) - \psi(r, \bar{x}) \bigg]\\
                                     & + e^{-\gamma j \left( \frac{\bar{x}}{N^{1/4} (\aa)}+kh \right)} \bigg[\psi\bigg(r, \bar{x} - j\frac{2}{N^{3/4}} \mathrm{sh}(\beta) \bigg) - \psi(r, \bar{x}) \bigg] \bigg\} \,,
\end{align*}
\end{allowdisplaybreaks}\\
where 
\begin{allowdisplaybreaks}
\begin{multline}\label{count'}
\vert A_N(i, j, k) \vert = 
\frac{N}{8}\bigg[1 + jk \mathrm{th}(\gamma h)+ik \mathrm{th}(\beta) \mathrm{th}(\gamma h)
+ ij \frac{\mathrm{th}(\beta) \mathrm{th}(\gamma h) \mathrm{sh}(\gamma h) + \mathrm{sh}(\beta)}{\aa}\\
+ k \frac{r}{N^{1/2}}
+ \frac{\bar{x}}{N^{1/4} (\aa)}\left(i
+ j \frac{\mathrm{ch}(\beta)}{ \mathrm{sh}(\beta)} + ijk \frac{\mathrm{sh}(\gamma h) (\bb)}{(\aa)^2}\right) \bigg]\,. 
\end{multline} 
\end{allowdisplaybreaks}\\
The procedure we applied to prove Corollary \ref{lemma Gspace scaling} leads us to the conclusion, once we leave all the terms coming from those processes we know collapsing in the infinite volume limit.
\fine

Now let consider the infinitesimal generator, $\mathcal{J}_N=N^{1/4}\mathcal{G}_N$, subject to the time-rescaling and apply it to the particular function $\psi (r_{N}(t), \widetilde{x}_{N}(t)) := \psi (r_{N}(t), \bar{x}_{N}(N^{1/4}t)) = \vert \widetilde{x}_{N}(t) \vert$.

The following decomposition holds
\[\begin{split}
\vert \widetilde{x}_{N}(t) \vert &= \vert \widetilde{x}_{N}(0) \vert + \int_0^t \mathcal{J}_N(\vert \widetilde{x}_{N}(s) \vert)ds + \mathcal{M}_{N, \vert \widetilde{x} \vert}^t\\
                           &\leq \vert \widetilde{x}_{N}(0) \vert + \int_0^t \vert \mathcal{J}_N(\vert \widetilde{x}_{N}(s) \vert) \vert ds + \mathcal{M}_{N, \vert \widetilde{x} \vert}^t\,,
\end{split}\]
with
\[
\mathcal{M}_{N, \vert \widetilde{x} \vert}^t = \int_0^t \sum_{i,j,k\in\mathscr{S}} \left\{ \overline{\nabla}^{(i)}[ \vert \widetilde{x}_{N}(s) \vert] \widetilde{\Lambda}^\sigma_N(i,j,k,ds) + \overline{\nabla}^{(j)}[ \vert \widetilde{x}_{N}(s) \vert]  \widetilde{\Lambda}^\omega_N(i,j,k,ds)  \right\}\,,
\]
where we have defined
\begin{align}\label{gradmodx'}
\overline{\nabla}^{(i)}[ \vert \widetilde{x}_{N}(t) \vert ] &:=  \bigg\vert \widetilde{x}_{N}(t) - i \frac{2}{N^{3/4}} \cosh(\gamma h) \bigg\vert - \vert \widetilde{x}_{N}(t) \vert \nonumber\\
&\\
\overline{\nabla}^{(j)}[ \vert \widetilde{x}_{N}(t) \vert ] &:=  \bigg\vert \widetilde{x}_{N}(t) - j \frac{2}{N^{3/4}} \sinh(\beta) \bigg\vert - \vert \widetilde{x}_{N}(t) \vert \nonumber
\end{align}
and 
\begin{align}\label{lambdatilde'}
\widetilde{\Lambda}^\sigma_N(i,j,k,dt) &:= \Lambda^\sigma_N(i,j,k,dt) - \underbrace{N^{1/4} \left\vert A(i,j,k,N^{1/4}t) \right\vert  e^{-\beta ij} dt}_{:= \lambda^\sigma (i,j,k,t)\,dt} \nonumber \\
&\\
\widetilde{\Lambda}^\omega_N(i,j,k,dt) \! &:= \! \Lambda^\omega_N(i,j,k,dt) \! - \! \underbrace{N^{1/4} \! \left\vert A(i,j,k,N^{1/4}t) \right\vert \!  e^{-\gamma j\left( \! \frac{\widetilde{x}_{N}(t)}{N^{1/4} [\cosh(\beta) + \cosh(\gamma h)]}+kh \!\right)} dt}_{:= \lambda^\omega (i,j,k,t)dt}\,. \nonumber
\end{align}
As we can clearly see, the quantities $\widetilde{\Lambda}^\cdot_N(i,j,k,dt)$ are the differences between the point processes $\Lambda^\cdot_N(i,j,k,dt)$, defined on \mbox{$\mathscr{S}^3 \times \mathbb{R}^+$}, and their intensities \mbox{$\lambda^\cdot (i,j,k,t)\,dt$}. \\
The counter $\left\vert A(i,j,k,N^{1/4}t) \right\vert$ is given in analogy with \eqref{count'}, replacing the variables $r$ and $\bar{x}$ with the stochastic processes $r_{N}(t)$ and $\widetilde{x}_{N}(t)$.\\
We recall that the expression of $\mathcal{G}_N$ is given by \eqref{IG4'}.\\
For $t \in [0, \tau_{N}^M]$ we can estimate $\vert \mathcal{J}_N(\,\vert \widetilde{x}_{N}(t) \vert\,) \vert $. We get
%\begin{allowdisplaybreaks}
\begin{multline*}
\vert \mathcal{J}_N(\,\vert \widetilde{x}_{N}(t) \vert\,) \vert =\bigg\vert 2N^{1/4} \mathrm{sgn}(\widetilde{x}_{N}(t)) \bigg\{ \frac{r_{N}(t)}{N^{1/4}} \mathrm{sh}(\beta) \mathrm{sh}(\gamma h)
 - \frac{\gamma^2 \vert \widetilde{x}_{N}(t) \vert^3}{2N^{1/2}} \frac{\mathrm{ch}(\beta) \mathrm{ch}(\gamma h)}{(\aa)^3}\\
 + \frac{\gamma^3 \vert \widetilde{x}_{N}(t) \vert^3}{6N^{1/2}} \frac{\mathrm{sh}(\beta)(\mathrm{ch}(\gamma h) - \mathrm{sh}(\gamma h) \mathrm{th}(\gamma h))}{ (\aa)^3}
 + \frac{\gamma^2 r_{N}(t) \vert \widetilde{x}_{N}(t) \vert^2}{2N^{3/4}} \frac{\mathrm{sh}(\beta) \mathrm{sh}(\gamma h)}{ (\aa)^2}\\
 + N^{1/4} R_s \left[\mathrm{sh}(\beta) \mathrm{ch}(\gamma h) - \mathrm{sh}(\beta) \mathrm{sh}(\gamma h) \mathrm{th}(\gamma h)\right]\\
 + N^{1/4} R_c \bigg[ \frac{r_{N}(t)}{N^{1/2}} \mathrm{sh}(\beta) \mathrm{sh}(\gamma h) - \frac{\vert \widetilde{x}_{N}(t) \vert}{N^{1/4}}  \frac{\mathrm{ch}(\beta)\mathrm{ch}(\gamma h)}{(\aa)} \bigg]  \bigg\} \bigg\vert %\\ \\
%\leq 2 \bigg\{M \mathrm{sh}(\beta) \mathrm{sh}(\gamma h) +  \frac{\gamma^2 M^3}{2N^{1/4}} \frac{\mathrm{ch}(\beta) \mathrm{ch}(\gamma h)}{ (\aa)^3}\\
% + \frac{\gamma^3 M^3}{6N^{1/4}} \frac{\mathrm{sh}(\beta)}{\mathrm{ch}(\gamma h) (\aa)^3}
% + \frac{\gamma^2 M^2}{2} \frac{\mathrm{sh}(\beta) \mathrm{sh}(\gamma h)}{(\cosh(\beta) + \cosh(\gamma h))^2}\\
% + \frac{\gamma^3 M^3}{6N^{1/4}} \mathrm{ch}\left[ \frac{\gamma M}{\aa} \right]\frac{ \mathrm{sh}(\beta)\left(\mathrm{ch}(\gamma h) +\mathrm{sh}(\gamma h) \mathrm{th}(\gamma h) \right)}{(\aa)^3}\\
% + \frac{\gamma^3 M^3}{6N^{1/4} } \mathrm{sh}\left[\frac{\gamma M}{\aa} \right]\frac{1}{(\aa)^3}\left(\mathrm{sh}(\beta) \mathrm{sh}(\gamma h)+\frac{M\mathrm{ch}(\beta) \mathrm{ch}(\gamma h)}{N^{1/4}(\aa)}\right)   \bigg\} \\
% \\
%\leq 2 \bigg\{M \mathrm{sh}(\beta) \mathrm{sh}(\gamma h) +  \frac{\gamma^2 M^3}{2} \frac{\mathrm{ch}(\beta) \mathrm{ch}(\gamma h)}{(\aa)^3}
% + \frac{\gamma^3 M^3}{6} \frac{\mathrm{sh}(\beta)}{\mathrm{ch}(\gamma h)(\aa)^3}\\
% + \frac{\gamma^2 M^2}{2} \frac{\mathrm{sh}(\beta) \mathrm{sh}(\gamma h)}{(\aa)^2}
% +\frac{\gamma^3 M^3}{6}\mathrm{ch}\left[\frac{\gamma M}{\aa} \right]\frac{ \mathrm{sh}(\beta)\left(\mathrm{ch}(\gamma h) +\mathrm{sh}(\gamma h) \mathrm{th}(\gamma h) \right)}{(\aa)^3}\\
% + \frac{\gamma^3 M^3}{6} \mathrm{sh}\left[ \frac{\gamma M}{\aa} \right]\frac{1}{(\aa)^3}\left(\mathrm{sh}(\beta) \mathrm{sh}(\gamma h)+\frac{M\mathrm{ch}(\beta) \mathrm{ch}(\gamma h)}{\aa}\right)   \bigg\} := C_{8}\,,
\end{multline*}
%\end{allowdisplaybreaks}
and, by using the Taylor expansions of the hyperbolic sine and cosine stopped at the second order (see \eqref{sinh}, \eqref{cosh}) and the estimates of their remainders (see \eqref{resto sinh}, \eqref{resto cosh}), it results to be bounded from above by a positive constant $C_8$, which is independent of $N$. Moreover,
%Since the following inclusions are valid,
\begin{allowdisplaybreaks}
\begin{align*}
\{&\tau_{N}^M \leq T\}\subseteq \bigg\{\sup_{0\leq t \leq T\wedge\tau_{N}^M}\{\vert \widetilde{x}_{N}(t) \vert, \vert \widetilde{y}_{N}(t) \vert, \vert \widetilde{z}_{N}(t) \vert, \vert \widetilde{u}_{N}(t) \vert, \vert \widetilde{v}_{N}(t) \vert, \vert \widetilde{w}_{N}(t) \vert \} \geq M\bigg\}\\
&\,\subseteq \bigg\{\sup_{0\leq t \leq T\wedge\tau_{N}^M} \vert \widetilde{x}_{N}(t) \vert \geq M\bigg\} \cup \bigg\{\sup_{0\leq t \leq T\wedge\tau_{N}^M} \vert \widetilde{y}_{N}(t) \vert \geq M\bigg\}
\cup \bigg\{\sup_{0\leq t \leq T\wedge\tau_{N}^M} \vert \widetilde{z}_{N}(t) \vert \geq M\bigg\}\\ 
&\,\quad\cup \bigg\{\sup_{0\leq t \leq T\wedge\tau_{N}^M} \vert \widetilde{u}_{N}(t) \vert \geq M\bigg\} 
\cup \bigg\{\sup_{0\leq t \leq T\wedge\tau_{N}^M} \vert \widetilde{v}_{N}(t) \vert \geq M\bigg\} \cup \bigg\{\sup_{0\leq t \leq T\wedge\tau_{N}^M} \vert \widetilde{w}_{N}(t) \vert \geq M\bigg\}\\
&\,\subseteq \bigg\{\sup_{0\leq t \leq T\wedge\tau_{N}^M} \vert \widetilde{y}_{N}(t) \vert \geq M\bigg\} \cup \bigg\{\sup_{0\leq t \leq T\wedge\tau_{N}^M} \vert \widetilde{z}_{N}(t) \vert \geq M\bigg\} \\
&\,\quad\cup \bigg\{\sup_{0\leq t \leq T\wedge\tau_{N}^M} \vert \widetilde{u}_{N}(t) \vert \geq M\bigg\} \cup  \bigg\{\sup_{0\leq t \leq T\wedge\tau_{N}^M} \vert \widetilde{v}_{N}(t) \vert \geq M\bigg\} 
\cup \bigg\{\sup_{0\leq t \leq T\wedge\tau_{N}^M} \vert \widetilde{w}_{N}(t) \vert \geq M\bigg\}\\ 
&\,\quad\cup \{ \vert \widetilde{x}_{N}(0) \vert \geq C_{9} \}
\cup\bigg[ \{ \vert \widetilde{x}_{N}(0) \vert \leq C_{9} \} \cap  \bigg\{\sup_{0\leq t \leq T\wedge\tau_{N}^M} \vert \widetilde{x}_{N}(t) \vert \geq C_{9} + TC_{8} + C_{10}\bigg\} \bigg] \\
&\,\subseteq \bigg\{\sup_{0\leq t \leq T\wedge\tau_{N}^M} \vert \widetilde{y}_{N}(t) \vert \geq M\bigg\} \cup \bigg\{\sup_{0\leq t \leq T\wedge\tau_{N}^M} \vert \widetilde{z}_{N}(t) \vert \geq M\bigg\}\\
&\,\quad\cup \bigg\{\sup_{0\leq t \leq T\wedge\tau_{N}^M} \vert \widetilde{u}_{N}(t) \vert \geq M\bigg\} \cup  \bigg\{\sup_{0\leq t \leq T\wedge\tau_{N}^M} \vert \widetilde{v}_{N}(t) \vert \geq M\bigg\} 
\cup \bigg\{\sup_{0\leq t \leq T\wedge\tau_{N}^M} \vert \widetilde{w}_{N}(t) \vert \geq M\bigg\}\\ 
&\,\quad\cup \{ \vert \widetilde{x}_{N}(0) \vert \geq C_{9} \} 
\cup \bigg\{\sup_{0\leq t \leq T\wedge\tau_{N}^M} \mathcal{M}_{N, \vert \widetilde{x} \vert}^t \geq C_{10}\bigg\}\,,
\end{align*}
\end{allowdisplaybreaks}\\
then we obtain the following inequality for the probability of the interested set 
\begin{multline*}
P\{\tau_{N}^M \leq T\} \leq P\bigg\{\sup_{0\leq t \leq T\wedge\tau_{N}^M} \vert \widetilde{y}_{N}(t) \vert \geq M\bigg\} + P\bigg\{\sup_{0\leq t \leq T\wedge\tau_{N}^M} \vert \widetilde{z}_{N}(t) \vert \geq M\bigg\}\\
 + P\bigg\{\sup_{0\leq t \leq T\wedge\tau_{N}^M} \vert \widetilde{u}_{N}(t) \vert \geq M\bigg\} +  P\bigg\{\sup_{0\leq t \leq T\wedge\tau_{N}^M} \vert \widetilde{v}_{N}(t) \vert \geq M\bigg\}
+ P\bigg\{\sup_{0\leq t \leq T\wedge\tau_{N}^M} \vert \widetilde{w}_{N}(t) \vert \geq M\bigg\}\\ 
+ P\{ \vert \widetilde{x}_{N}(0) \vert \geq C_{9} \}
 + P\bigg\{\sup_{0\leq t \leq T\wedge\tau_{N}^M} \mathcal{M}_{N, \vert \widetilde{x} \vert}^t \geq C_{10}\bigg\}\,.
\end{multline*}
We estimate the seven terms of the right-hand side of the inequality. 
\begin{itemize}
\item[$\RHD$] for any $\varepsilon >0$, thanks to the fact that the process $\widetilde{z}_{N}(t)$ collapses we have 
\[P\bigg\{\sup_{0\leq t \leq T\wedge\tau_{N}^M} \vert \widetilde{z}_{N}(t) \vert \geq M\bigg\}\leq\varepsilon\,,\]
where we set $M := C_7 \left( N^{1/(8d)-1/16} \vee  N^{-1/8} \right)$ (see \eqref{thcorclps}) and analogous relations hold for all the other processes $\widetilde{y}_{N}(t)$, $\widetilde{u}_{N}(t)$, $\widetilde{v}_{N}(t)$, $\widetilde{w}_{N}(t)$, with proper constants;
\item[$\RHD$] from \eqref{crtdir'} we get 
$$E[\widetilde{x}_{N}(0)]=N^{1/4}E \Big[ \cosh(\gamma h) m_{\rho_N(0)}^{\underline{\sigma}} + \sinh(\beta) m_{\rho_N(0)}^{\underline{\omega}} \Big]\,,$$
which is a linear combination of sample averages. Since at time $t=0$ the spins are distributed according to a product measure and such that $m_0^{\sigma}=m_0^{\omega}=0$, thanks to the Central Limit Theorem we can conclude that
\[ E[\vert \widetilde{x}_{N}(0) \vert] \leq \left[ \cosh(\gamma h) \sqrt{\mathrm{Var}(\sigma_1(0))} + \sinh(\beta) \sqrt{\mathrm{Var}(\omega_1(0))} \right]N^{-1/4}\,,\]
and so, in the limit as $N \longrightarrow +\infty$, we have convergence to zero in $L^1$ and then in probability. Therefore
\[ P\{ \vert \widetilde{x}_{N}(0) \vert \geq C_{9} \}\leq \varepsilon\,, \]
for any $\varepsilon >0$, for every $N$ and for a sufficiently large $C_{9}$; 
\item[$\RHD$] we reduce to deal with $E[(\mathcal{M}_{N, \vert \widetilde{x} \vert}^T)^2]$; in fact, Doob's ``maximal inequality in $L^p$'' (case $p=2$) for martingales (we refer to Chapter VII, Section 3 of \cite{Shi96}) tells us that $$P\,\bigg\{\sup_{0\leq t \leq T\wedge\tau_{N}^M}\mathcal{M}_{N, \vert \widetilde{x} \vert}^t \geq C_{10}\bigg\} \leq \frac{E\left[\left(\mathcal{M}_{N, \vert \widetilde{x} \vert}^T\right)^2\right]}{(C_{10})^2}\, .$$
Hence, remembering \eqref{lambdatilde'} and \eqref{gradmodx'}, we are able to compute
\begin{allowdisplaybreaks}
\begin{align*}
E\left[\left(\mathcal{M}_{N, \vert \widetilde{x} \vert}^T\right)^2\right] &= E\bigg[\int_0^T \sum_{i,j,k\in\mathscr{S}} \Big\{ \Big[ \overline{\nabla}^{(i)}[\vert \widetilde{x}_{N}(t) \vert]\Big]^2 \lambda^\sigma (i,j,k,t) dt\\
& \hspace{115pt} + \Big[ \overline{\nabla}^{(j)}[\vert \widetilde{x}_{N}(t) \vert]\Big]^2 \lambda^\omega (i,j,k,t) dt \Big\} \bigg]\\
& \leq E\bigg[\int_0^T \frac{4}{N^{3/2}} \cosh^2(\gamma h) N^{1/4} \sup_{i,j,k \in \mathscr{S}} \vert A(i,j,k, N^{1/4}t) \vert \,e^{\beta} dt\\
& \qquad\qquad + \frac{4}{N^{3/2}} \sinh^2(\beta) N^{1/4} \sup_{i,j,k \in \mathscr{S}} \vert A(i,j,k, N^{1/4}t) \vert \,e^{\gamma(1+h)}dt\bigg]\\
& \leq E\bigg[\int_0^T \frac{4}{N^{5/4}} N \left( \cosh^2(\gamma h) e^{\beta} + \sinh^2(\beta) e^{\gamma(1+h)} \right) dt\bigg]\\
&\leq 4T \left( \cosh^2(\gamma h) e^{\beta} + \sinh^2(\beta) e^{\gamma(1+h)} \right) =: C_{11}  \,,
\end{align*}
\end{allowdisplaybreaks}\\
with $C_{11}$ positive constant independent on $N$ and $M$.\\ 
We have established that, if we choose $C_{10}\geq \sqrt{C_{11}/\varepsilon}$, then 
\[ P\bigg\{\sup_{0\leq t \leq T\wedge\tau_{N}^M} \mathcal{M}_{N, \vert \widetilde{x} \vert}^t \geq C_{10}\bigg\} \leq \varepsilon\,. \] 
\end{itemize}

In summary, we proved the inequality we were looking for; in fact
\[  P\left\{\tau_{N}^M\leq T \right\} \leq 7\varepsilon := \epsilon\,. \] 
%%%%%%%%%%%%%%%%%%%%%%%%%%%%%%%%%%%%%%%%%%%%%%%%%%%%%%%%%%%%%%%%%%%%%%%%%%%%%%%
This completes the first part of the proof.

Now, we are going to show that in the limit of infinite volume, when $t \in [0,T]$, the process $\widetilde{x}_{N}(t)$ admits a limiting process and we are going to compute it.\\

%%%%%%%%%%%%%%%%%%%%%%STEP 4%%%%%%%%%%%%%%%%%%%%%%%
\textbf{STEP 4.} First, we need to prove the tightness of the sequence $\{\widetilde{x}_{N}(t)\}_{N \geq 1}$. This property implies the existence of convergent subsequences. Secondly, in the last step, we will verify that all the convergent subsequences have the same limit and hence also the sequence $\{\widetilde{x}_{N}(t)\}_{N \geq 1}$ must converge to that limit. 

\begin{LEMMA}\label{LmmTight'}
The sequence $\{\widetilde{x}_{N}(t)\}_{N \geq 1}$ is tight.
\end{LEMMA}

\Proof
In the case we are working with processes with laws on $\mathcal{D}[0,T]$, we can give a characterization of the tightness in terms of those processes (through their distributions). In fact, as we can read in \cite{CoEi88}, we have: \medskip\\
``A sequence of processes $\{\widetilde{x}_N(t)\}_{N \geq 1}$ with laws $\{\mathcal{P}_N\}_{N \geq 1}$ on $\mathcal{D}[0,T]$ is tight if:
\begin{enumerate}
\item for every $\varepsilon>0$ there exists $M>0$ such that
\begin{equation}\label{T1}
\sup_{N}P\bigg\{\sup_{t\in[0,T]} \vert \widetilde{x}_N(t)\vert \geq M\bigg\} \leq \varepsilon\,,
\end{equation}
\item for every $\varepsilon>0$ and $\alpha >0$ there exists $\delta > 0$ such that
\begin{equation}\label{T2}
\sup_{N}\sup_{0\leq \tau_1 \leq \tau_2 \leq (\tau_1 + \delta) \wedge T }P\{\vert \widetilde{x}_N(\tau_2) - \widetilde{x}_N(\tau_1) \vert \geq \alpha\} \leq \varepsilon\,,
\end{equation}
where $\tau_1$ and $\tau_2$ are stopping times adapted to the filtration generated by the process $\widetilde{x}_N$.''
\end{enumerate}

We must verify the conditions \eqref{T1} and \eqref{T2} hold. Since we have already shown that, for every $\epsilon>0$ the inequality $P \left\{\tau_N^M \leq T \right\} \leq \epsilon$ holds for $M$ sufficiently large and uniformly in $N$, it is enough to show tightness for the stopped process
\[
\left\{\widetilde{x}_N \left( t \wedge \tau_N^M \right) \right\}_{N \geq 1} \,.
\]
We showed before the validity of the following inclusion
\[\bigg\{\sup_{0\leq t \leq T\wedge\tau_{N}^M} \vert \widetilde{x}_{N}(t) \vert \geq M\bigg\} \subseteq \{ \vert \widetilde{x}_{N}(0) \vert \geq C_{9} \} \cup \bigg\{\sup_{0\leq t \leq T\wedge\tau_{N}^M} \mathcal{M}_{N, \vert \widetilde{x} \vert}^t \geq C_{10}\bigg\}\,,\]
therefore
 \[ \sup_{N} P\bigg\{\sup_{0\leq t \leq T\wedge\tau_{N}^M} \vert \widetilde{x}_{N}(t) \vert \geq M\bigg\} \leq 2\varepsilon \] 
and so we obtain \eqref{T1}. 

Now let us deal with \eqref{T2}. We notice that
\[ \vert  \widetilde{x}_{N}(\tau_2) - \widetilde{x}_{N}(\tau_1) \vert = \bigg\vert \int_{\tau_1}^{\tau_2} \mathcal{J}_N (\widetilde{x}_{N}(u))du + \mathcal{M}_{N, \vert \widetilde{x} \vert}^{\tau_1,\tau_2} \bigg\vert \,,\]
where we have denoted
\[
\mathcal{M}_{N, \vert \widetilde{x} \vert}^{\tau_1,\tau_2} = -\,\frac{2}{N^{3/4}}\int_{\tau_1}^{\tau_2} \sum_{i,j,k\in\mathscr{S}} \left( i \cosh(\gamma h) \widetilde{\Lambda}_N^\sigma(i,j,k,du) + j \sinh(\beta) \widetilde{\Lambda}_N^\omega(i,j,k,du) \right) 
\]
as in definition \eqref{lambdatilde'}. Thus, 
\[\{\vert \widetilde{x}_{N}(\tau_2) - \widetilde{x}_{N}(\tau_1) \vert \geq \alpha\} \subseteq \bigg\{\underbrace{\bigg\vert \int_{\tau_1}^{\tau_2} \mathcal{J}_N \vert \widetilde{x}_{N}(u)
\vert du \bigg\vert}_{\leq C_{8}(\tau_2 - \tau_1)} + \left\vert \mathcal{M}_{N, \vert \widetilde{x} \vert}^{\tau_1,\tau_2} \right\vert \geq \alpha \bigg\} \subseteq \left\{ \left\vert \mathcal{M}_{N, \vert \widetilde{x} \vert}^{\tau_1,\tau_2} \right\vert \geq \overline{C}_{10} \right\} \]
and then, applying Chebyscev Inequality to the last right-hand side of the previous inclusions, we get
\[
\sup_{0\leq \tau_1 \leq \tau_2 \leq (\tau_1 + \delta) \wedge T }P \left\{\left\vert \mathcal{M}_{N, \vert \widetilde{x} \vert}^{\tau_1,\tau_2} \right\vert \geq \overline{C}_{10} \right\} \leq (\overline{C}_{10})^{-2} \sup_{0\leq \tau_1 \leq \tau_2 \leq (\tau_1 + \delta) \wedge T } E \left[ \left(\mathcal{M}_{N, \vert \widetilde{x} \vert}^{\tau_1,\tau_2} \right)^2 \right]
\]
Observing that $\mathcal{M}_{N, \vert \widetilde{x} \vert}^{t}$ is a zero mean martingale, by Doob's Optional Sampling theorem, we obtain
\begin{align*}
(\overline{C}_{10})^{-2}  & \sup_{0\leq \tau_1 \leq \tau_2 \leq (\tau_1 + \delta) \wedge T } E \left[ \left(\mathcal{M}_{N, \vert \widetilde{x} \vert}^{\tau_1,\tau_2} \right)^2 \right]\\
& =(\overline{C}_{10})^{-2} \!\! \sup_{0\leq \tau_1 \leq \tau_2 \leq (\tau_1 + \delta) \wedge T } E \left[ \left(\mathcal{M}_{N, \vert \widetilde{x} \vert}^{\tau_2} \right)^2 \!\! - \left(\mathcal{M}_{N, \vert \widetilde{x} \vert}^{\tau_1} \right)^2 \right] \\
& \leq (\overline{C}_{10})^{-2}\sup_{0\leq \tau_1 \leq \tau_2 \leq (\tau_1 + \delta) \wedge T } 4 (\tau_2 - \tau_1)\Big[ \cosh^2(\gamma h) e^{\beta}%\\
%& \hspace{135pt} 
+ \sinh^2(\beta) e^{\gamma(1+h)} \Big]\\
%& 
& \leq (\overline{C}_{10})^{-2} \underbrace{4 \Big[ \cosh^2(\gamma h) e^{\beta} + \sinh^2(\beta) e^{\gamma(1+h)} \Big]}_{:=\overline{C}_{11}}\delta \,.
%\end{split}\]
\end{align*}
Finally, we can conclude that
\[\begin{split}
\sup_{N}\sup_{0\leq \tau_1 \leq \tau_2 \leq (\tau_1 + \delta) \wedge T }P\{\vert \widetilde{x}_{N}(\tau_2) - \widetilde{x}_{N}(\tau_1) \vert \geq \alpha\}
&\leq \sup_{N}\sup_{0\leq \tau_1 \leq \tau_2 \leq (\tau_1 + \delta) \wedge T }P \left\{ \left\vert \mathcal{M}_{N, \vert \widetilde{x} \vert}^{\tau_1,\tau_2} \right\vert \geq \overline{C}_{10} \right\}\\
&\leq (\overline{C}_{10})^{-2} \overline{C}_{11} \, \delta = O(\delta)
\end{split}\]
and the proof is complete.
\fine

%%%%%%%%%%%%%%%%%%%%%%%%%%%%%%%%%STEP 5%%%%%%%%%%%%%%%%%%%%%%%%%%%%%%%%%
\textbf{STEP 5.} We prove now that all the convergent subsequences have the same limit and so the sequence itself converges to that limit and this concludes the proof of the Theorem.\\
With abuse of notation, let $\{\widetilde{x}_{n}(t)\}_{n \geq 1}$ denote one of such a subsequence and let $\psi \in \mathcal{C}^3_b$ be a function of the type $\psi(r_n(t), \widetilde{x}_n(t)) = \psi (\widetilde{x}_n(t))$. The following decomposition holds
\begin{equation}\label{MartPbFinale'}
\psi(\widetilde{x}_{n}(t)) - \psi(\widetilde{x}_{n}(0)) = \int_0^t \mathcal{J}_n \psi(\widetilde{x}_{n}(u))du + \mathcal{M}_{n,\psi}^t\,,
\end{equation}
where
\begin{multline*}
\mathcal{J}_n \psi(\widetilde{x}_{n}(t))  =%& 
\,2 \psi_{\widetilde{x}}\,\bigg[r_{n}(t) \sinh(\beta) \sinh(\gamma h) - \frac{\gamma^2 (\widetilde{x}_{n}(t))^3}{2} \frac{\cosh(\beta) \cosh(\gamma h)}{n^{1/4} [\cosh(\beta) + \cosh(\gamma h)]^3}\\
                                 %& 
                                 \quad\qquad 
                                 + \frac{\gamma^3 (\widetilde{x}_{n}(t))^3}{6} \frac{\sinh(\beta)[\cosh(\gamma h) - \sinh(\gamma h) \tanh(\gamma h)]}{n^{1/4} [\cosh(\beta) + \cosh(\gamma h)]^3}\\
                                 %& \qquad 
                                 + \frac{\gamma^2 r_n(t) (\widetilde{x}_{n}(t))^2}{2} \frac{\sinh(\beta) \sinh(\gamma h)}{h n^{3/4} [\cosh(\beta) + \cosh(\gamma h)]^2} \bigg]  + o_M(1)\,,
\end{multline*}
which, as usual, is $\mathcal{G}_N$ (see \eqref{IG4'}) rescaled by a power $N^{1/4}$ and applied to the particular  function $\psi(r_{n}(t), \widetilde{x}_{n}(t))=\psi(\widetilde{x}_{n}(t))$. The remainder $o_M(1)$ goes to zero as $n \longrightarrow +\infty$, uniformly in $M$. \\
If we compute the limit as $n\longrightarrow +\infty$, remembering that a Central Limit Theorem applies to $r_{n}(t)$, we have 
\[ \mathcal{J}_n \psi(\widetilde{x}_{n}(t)) \xrightarrow[\quad w \quad]{n\rightarrow +\infty} \mathcal{J} \psi(\widetilde{x}(t))\,, \]
with
\[ \mathcal{J} \psi(\widetilde{x}(t)) = 2 \, \mathscr{H} \, \sinh(\beta) \sinh(\gamma h)\; \psi_{\widetilde{x}}  \]
and $\mathscr{H}$ is a Standard Gaussian random variable. Then, because of \eqref{MartPbFinale'}, we obtain 
\[  \mathcal{M}_{n,\psi}^t \xrightarrow[\quad w \quad]{n \rightarrow +\infty}  \mathcal{M}_{\psi}^t := \psi(\widetilde{x}(t)) - \psi(\widetilde{x}(0)) - \int_0^t \mathcal{J} \psi(\widetilde{x}(u))du\,. \] 
We must prove the following Lemma

\begin{LEMMA}\label{LmmMart}
$M_\psi^t$ is a martingale (with respect to $t$); in other words, for all $s,t \in [0,T]$, $s\leq t$ and for all measurable and bounded functions  $g(\widetilde{x}([0,s]))$ the following identity holds: 
\begin{equation}\label{mart}
E[ \mathcal{M}_{\psi}^t g(\widetilde{x}([0,s]))] = E[ \mathcal{M}_{\psi}^s g(\widetilde{x}([0,s]))]\,.
\end{equation}
\end{LEMMA}

\Proof
It is sufficient to show that $ \{ \mathcal{M}_{n,\psi}^t \}_{n \geq 1}$ is a uniformly integrable sequence of random variables. Let us suppose we have already proved this property holds and see that \eqref{mart} is satisfied.\\
Since $\mathcal{M}_{n,\psi}^t$ is a martingale (with respect to $t$) for every $n$, we have that, for all $s,t \in [0,T]$, $s\leq t$ and for all measurable and bounded functions $g(\widetilde{x}([0,s]))$,
\[
E[ \mathcal{M}_{n,\psi}^t g(\widetilde{x}([0,s]))] = E[ \mathcal{M}_{n,\psi}^s g(\widetilde{x}([0,s]))]
\]  
and then
\[
\lim_{n \rightarrow +\infty} E[ \mathcal{M}_{n,\psi}^t g(\widetilde{x}([0,s]))] = \lim_{n \rightarrow +\infty} E[ \mathcal{M}_{n,\psi}^s g(\widetilde{x}([0,s]))]\,.
\]
But $\{ \mathcal{M}_{n,\psi}^t \}_{n \geq 1}$ is a sequence of uniformly integrable random variables, hence it converges in $L^1$ (for instance, see ~\cite{Shi96}). Moreover, we know the distribution of its $L^1$-limit, since we already know its weak-limit. Thus,
\begin{align*}
E[ \mathcal{M}_{\psi}^t g(\widetilde{x}([0,s]))] &= E \left[ \lim_{n \rightarrow +\infty} \mathcal{M}_{n,\psi}^t g(\widetilde{x}([0,s])) \right] = \lim_{n \rightarrow +\infty} E[ \mathcal{M}_{n,\psi}^t g(\widetilde{x}([0,s]))]\\
&= \lim_{n \rightarrow +\infty} E[ \mathcal{M}_{n,\psi}^s g(\widetilde{x}([0,s]))] =  E \left[ \lim_{n \rightarrow +\infty} \mathcal{M}_{n,\psi}^s g(\widetilde{x}([0,s])) \right]\\
&= E[ \mathcal{M}_{\psi}^s g(\widetilde{x}([0,s]))]
\end{align*}
and the conclusion follows.

It remains to check  that $\{ \mathcal{M}_{n,\psi}^t \}_{n \geq 1}$ is a uniformly integrable family. A sufficient condition for the uniform integrability is the existence of $p>1$ such that $\sup_n E[\vert  \mathcal{M}_{n,\psi}^t \vert ^p]< +\infty$ (see again ~\cite{Shi96}). \\
If we define 
\begin{align*}
\overline{\nabla}^{(i)}[\psi(\widetilde{x}_{n}(t))] &:=  \psi \bigg( \widetilde{x}_{n}(t) - i \frac{2}{n^{3/4}} \cosh(\gamma h) \bigg) - \psi (\widetilde{x}_{n}(t)) \\
\overline{\nabla}^{(j)}[\psi(\widetilde{x}_{n}(t))] &:=  \psi \bigg( \widetilde{x}_{n}(t) - j \frac{2}{n^{3/4}} \sinh(\beta) \bigg) - \psi (\widetilde{x}_{n}(t))\,,
\end{align*}
it yields
\begin{allowdisplaybreaks}
\begin{align*}
&E \left[ \left( \mathcal{M}_{n,\psi}^t \right)^2 \right]=\\ 
&\quad=E\bigg[ \int_0^t \sum_{i,j,k\in\mathscr{S}} \Big\{ \Big[ \overline{\nabla}^{(i)}[\psi(\widetilde{x}_{n}(s))] \Big]^2 \lambda^\sigma (i,j,k,s)ds
%\\ &\hspace{130pt} 
+ \Big[ \overline{\nabla}^{(j)}[\psi(\widetilde{x}_{n}(s))] \Big]^2 \lambda^\omega (i,j,k,s)ds \Big\} \bigg]
\\ 
& 
\quad\leq n^{5/4} E\bigg[ \int_0^t \sum_{i,j \in \mathscr{S}} \Big\{ \Big[ \overline{\nabla}^{(i)}[\psi(\widetilde{x}_{n}(s))] \Big]^2 e^{\beta}%\\
%&\hspace{150pt} 
+ \Big[ \overline{\nabla}^{(j)}[\psi(\widetilde{x}_{n}(s))] \Big]^2 e^{\gamma (1+h)} \Big\} ds \bigg] =(\ast)\\
\intertext{we expand the function $\psi$ around $\widetilde{x}_{n}(t)$ with the Taylor expansion stopped at the first order and with remainder $R$, $\overline{R}$ such that 
\begin{align*}
\vert R \vert &
\leq \frac{1}{2} \sup \left\{ \vert \psi_{\widetilde{x} \widetilde{x}}(z) \vert : z \in \left[ \widetilde{x}_{n}(t), \widetilde{x}_{n}(t) - i \frac{2}{n^{3/4}} \cosh(\gamma h) \right] \right\} \frac{4}{n^{3/2}} \cosh^2(\gamma h)\\
\vert \overline{R} \vert &
\leq \frac{1}{2} \sup \left\{ \vert \psi_{\widetilde{x} \widetilde{x}}(z) \vert : z \in \left[ \widetilde{x}_{n}(t), \widetilde{x}_{n}(t) - j \frac{2}{n^{3/4}} \sinh(\beta) \right] \right\} \frac{4}{n^{3/2}} \sinh^2(\beta)
\end{align*}
and, moreover, we recall that $\psi \in \mathcal{C}^3_b$ and so $\vert \psi_{\widetilde{x}} \vert \leq k_1$ and $\vert \psi_{\widetilde{x} \widetilde{x}}\vert \leq k_2$; therefore,}
& 
(\ast)= n^{5/4}\, E\bigg[ \int_0^t \bigg\{ \sum_{i \in \mathscr{S}} \bigg[ - i \frac{2}{n^{3/4}}\cosh(\gamma h)\psi_{\widetilde{x}} + R \,\bigg]^2 e^{\beta}\\
& %\hspace{95pt} 
\qquad\qquad\qquad\quad + \sum_{j \in \mathscr{S}} \bigg[ - j \frac{2}{n^{3/4}}\sinh(\beta)\psi_{\widetilde{x}} + \overline{R} \bigg]^2 e^{\gamma (1+h)} \bigg\} ds \,\bigg]\\
&
\leq n^{5/4}  \,E\bigg[ e^{\beta} \!\! \int_0^t \! \sup_{i \in \mathscr{S}} \bigg( \frac{4}{n^{3/2}} \cosh^2(\gamma h) \psi_{\widetilde{x}}^2 - i \frac{4}{n^{3/4}} \cosh(\gamma h) \psi_{\widetilde{x}} R +  R^2 \! \bigg) ds\\
& \qquad\qquad\quad 
+ e^{\gamma (1+h)} \int_0^t \sup_{j \in \mathscr{S}} \bigg( \frac{4}{n^{3/2}} \sinh^2(\beta) \psi_{\widetilde{x}}^2 - \frac{4}{n^{3/4}} \sinh(\beta) \psi_{\widetilde{x}} \overline{R} + \overline{R}^2 \bigg) ds\, \bigg]\\
&
\leq n^{5/4}\, E\bigg[ e^{\beta} \int_0^t \bigg( \frac{4}{n^{3/2}} k_1^2 \cosh^2(\gamma h) + \frac{8}{n^{9/4}} k_1 k_2 \cosh^3(\gamma h)%\\ 
%& \hspace{50pt} 
+ \frac{4}{n^{3}} k_2^2 \cosh^4(\gamma h) \bigg) ds\\ 
& \qquad\qquad\quad + e^{\gamma (1+h)} \int_0^t \bigg( \frac{4}{n^{3/2}} k_1^2 \sinh^2(\beta)%\\
%& \hspace{113pt} 
+ \frac{8}{n^{9/4}} k_1 k_2 \sinh^3(\beta) + \frac{4}{n^{3}} k_2^2 \sinh^4(\beta) \bigg) ds \,\bigg] \\
&
\leq 4 T \left[ e^{\beta} \cosh^2(\gamma h) (k_1 + \cosh(\gamma h) k_2)^2%\\
%& \hspace{147pt} 
+ e^{\gamma (1+h)} \sinh^2(\beta) (k_1 + \sinh(\beta) k_2)^2 \right] \,,
\end{align*}
\end{allowdisplaybreaks}\\
since $t<T$; then $ \mathcal{M}_{n,\psi}^t$ is uniformly integrable. 
\fine

Now, the proof is easy  to complete. $ \mathcal{M}_{n,\psi}^t$ solves the martingale problem with infinitesimal generator $\mathcal{J}$, admitting a unique solution, and hence we have shown that all the convergent subsequences have the same limit and so the sequence itself converges to that limit.

\appendix

%\bibliographystyle{plain}
%\bibliography{Bibliography}

\end{document}